\documentclass{svmult}
\usepackage{makeidx}         
\usepackage{graphicx}        
\usepackage{multicol}        
\usepackage[bottom]{footmisc}

\makeindex             

\renewcommand{\theequation}{\arabic{section}.\arabic{equation}}
\newtheorem{thm}{Theorem}[section]
\newtheorem{lem}{Lemma}[section]

\newtheorem{defi}{Definition}[section]

\begin{document}

\baselineskip15pt

\title*{Random Periodic Solutions
 of SPDEs via Integral Equations and Wiener-Sobolev Compact Embedding }
\titlerunning{Random Periodic Solutions}
\author{Chunrong Feng and Huaizhong Zhao}
\authorrunning{C. R. Feng and H. Z. Zhao}
\institute{ Department of Mathematical Sciences, Loughborough
University, LE11 3TU, UK\\
 \texttt{C.Feng@lboro.ac.uk}, \texttt{H.Zhao@lboro.ac.uk}}

\maketitle
\newcounter{bean}
\begin{abstract}
In this paper, we study the existence of random periodic solutions for
semilinear SPDEs on a bounded domain with a smooth boundary. We identify them
as the solutions of coupled forward-backward infinite horizon stochastic
integral equations on $L^2(D)$ in general cases. For this we use Mercer's Theorem and eigenvalues and eigenfunctions of the second order differential operators in the infinite horizon integral equations. We then use the argument of the relative compactness of Wiener-Sobolev spaces in 
$C^0([0, T], L^2(\Omega\times D))$ and generalized Schauder's fixed point theorem
to prove the existence of a solution of the integral equations.
This is the first paper in literature to study random periodic solutions of SPDEs. Our result is also new in finding semi-stable stationary solution for non-dissipative SPDEs, while in literature the classical method is to use the pull-back technique so researchers were only able to find stable stationary solutions for dissipative systems.\\

%

\vskip25pt
\end{abstract}

 \renewcommand{\theequation}{\arabic{section}.\arabic{equation}}

\section{Introduction}

Dynamics of nonlinear differential equations, both deterministic and stochastic, are complex. It is of great importance to understand these complexities. Mathematicians have made enormous progress in understanding these complexities for deterministic systems, both of finite dimensional and infinite dimensional. Understanding the complexities of  stochastic systems are far from clear even for stationary solutions. The concept of  stationary solutions is the stochastic counter part of fixed points to deterministic dynamical systems. A fixed point is the simplest equilibrium and large time limiting set of a deterministic dynamical system. A periodic solution is a more complicated limiting set. The theory of periodic solutions has played a central role in the study of the complex behaviour of a dynamical system. They are relatively simple trajectories themselves. However, their existence and construction is a challenging problem in the study of dynamical systems. The study has occupied a central role in the theory of dynamical system since the seminal work Henri Poincar\'e \cite{poincare}. Periodic solutions of partial differential equations of parabolic type has been studied by a number of authors, Vejvoda \cite{vejvoda}, Fife \cite{fife}, Hess \cite{hess}, Lieberman \cite {lieberman1}, \cite{lieberman2}, to name but a few. From periodic solutions, more complicated solutions can be built in. Since the theory of the existence of the solution of the stochastic differential equations (SDEs) and stochastic partial differential equations (SPDEs) become better understood (Da Prato and Zabczyk \cite{da-za1}, Pr\'ev$\hat{\rm {o}}$t and R$\ddot{\rm {o}}$ckner \cite{prevot-rockner}) we need to study more detailed question about the behaviour of solutions of SDEs and SPDEs. Mathematicians have been very much interested in the study of the existence of stationary solutions of SDEs and SPDEs, and invariant manifolds near stationary solutions. For results about SPDEs, see Sinai \cite {si1}, \cite{ si2}, Mattingly \cite{mattingly}, E, Khanin, Mazel and Sinai \cite{kh-ma-si}, Caraballo, Kloeden and Schmalfuss \cite{kloeden}, Liu and Zhao \cite{li-zh}, Zhang and Zhao \cite{zh-zh}, \cite{zh-zh1}, Duan, Lu and Schmalfuss \cite{du-lu-sc1}, \cite{du-lu-sc2}, Mohammed, Zhang and Zhao \cite{mo-zh-zh}, Lian and Lu \cite{lian-lu},  though there are still many problems that need to be understood. 
In 
literature, there were only few works on periodicity of stochastic systems. For linear stochastic differential 
equations with periodic coefficients in the sense of distribution, see Chojnowska-Michalik \cite {cho1}, \cite{cho2}, and for one-dimensional 
random mappings, see Kl\"unger \cite{klunger}.
We began to address the problem of pathwise random periodic solutions to SDEs in Zhao and Zheng \cite {zh-zheng}, Feng, Zhao and Zhou \cite{feng-zhao-zhou}. In this context, first we would like to motivate the reader with the following question. Consider a deterministic evolution equation on a Hilbert space $H$,
\begin{eqnarray}\label{1.1}
{du\over dt}=Au+f(u).
\end{eqnarray}
Assume it has a periodic solution of periodic $\tau$, $Z: (-\infty,\infty)\to H$ such that $Z(t+\tau)=Z(t)$, for any $t\in(-\infty, \infty)$. Now we consider the following stochastic differential equation, which can be regarded formally as the random perturbation of (\ref {1.1}) with a white noise perturbation:
\begin{eqnarray}\label{1.2}
du=(Au+f(u))dt +g(u)dW(t).
\end{eqnarray}
Here $W$ is a two-sided Brownian motion on a probability space $(\Omega, {\cal F}, P)$ valued in a Hilbert space $K$ and $g: H\to L_2(K, H)$ taking values in the space of 
Hilbert-Schmidt operators. Assume the solution of such an equation with a given initial condition exists and is unique. Such an equation has been considered in literature for many 
SDEs and SPDEs. The question to ask is: does equation (\ref{1.2}) still possess a periodic solution? Of course the answer is definitely no in general if we think periodic solution a 
close trajectory as in the deterministic sense. But a close trajectory is not the right notion of random periodic solution  to stochastic systems, just like the deterministic fixed point is 
not a right notion for stochastic systems. One can not expect that, in general, equation (\ref{1.2}) has a solution such that $u(t+\tau)=u(t)$ unless in a very special situation. There is an interaction between the periodic solution and the noise. Intuitively, the periodic solution has tendency to make trajectories of the random dynamical system following a periodic circle, at least in the dissipative case. The noise tends to make trajectories spreading out. Understanding of this kind of phenomenon was attempted by considering first linear approximation in physics literature, assuming the deterministic macroscopic equation has a periodic solution (see e.g. \cite{physical}). 
Note the following observation: let
\begin{eqnarray*}
u(t)=Z(t)+v(t).
\end{eqnarray*}
Then $v$ satisfies 
\begin{eqnarray}\label{1.3}
dv(t)=(Av(t)+b(t,v(t)))dt+\sigma(t, v(t))dW(t),
\end{eqnarray}
where 
\begin{eqnarray*}
b(t,v)&=&f(Z(t)+v)-f(Z(t)),\\
\sigma(t,v)&=&g(Z(t)+v).
\end{eqnarray*}
Note $b,\sigma$ are periodic function in $t$, i.e. $b(t+\tau,v)=b(t,v)$ and $\sigma(t+\tau,v)=\sigma(t,v)$ for any $t\in R$ and $v\in H$. Now the question is reduced to the study of 
the random periodic solution of equation (\ref{1.3}) with periodic coefficients. In fact, this kind of stochastic differential equations with periodic coefficients arises in modelling many physical problems. For example, it was considered in climate dynamics literature that mid-latitude oceans can be modelled by time periodic wind forcing when one takes into account the seasonal cycles in winds. But a more realistic model should include a stochastic effects (\cite{chekroun}). The periodic solution is naturally extended to the notion of the random periodic solution to equation 
such as equation (\ref{1.3}) with periodic coefficients by \cite{feng-zhao-zhou}.
If the periodic solution $Z$ of Equation (\ref {1.1}) is exponentially stable and the noise is reasonably small in Equation (\ref{1.3}) ($g(u)$ is Lipschitz in $u$ and the Lipschitz constant is reasonably small), we can construct a stable random periodic solution to equation (\ref{1.3}) therefore obtain a random periodic solution of equation (\ref{1.2}). But in the non-dissipative case  that equation (\ref{1.1}) has a periodic solution $Z$ of period $\tau$, not stable but semi-stable, the situation is more complicated. Pull-back 
and Poincar\'e mapping approaches do not seem working easily in this situation.

In {\cite{feng-zhao-zhou}, we proved in the case that $H=R^d$ and $A$ is hyperbolic the existence of random periodic solution of Equation (\ref{1.3}) is equivalent to the existence of a solution of an infinite horizon $(-\infty,\infty)$ integral equation. In fact, the result holds in both finite and infinite dimensional spaces, though we only gave the proof in the $R^d$ case. Furthermore, we extended the Schauder fixed point theorem to the case when the subspace of the Banach space is not closed and the Wiener-Sobolev compactness theorem to the relative compactness on the space $C([0,T], L^2(dP))$. Then we proved the existence of a solution of the infinite horizon integral equation.

In this paper, we continue to push this new idea to the following stochastic partial differential equation of 
parabolic type on a bounded domain $D\subset R^d$ with a smooth boundary:
\begin{eqnarray}\label{march17a}
du(t,x)&=&{\cal L} u(t,x)\,dt+F(t,u(t,x))\, dt+\sum_{k=1}^\infty \sigma_k(t)\phi_k(x)dW^k(t), \ \ \ \ t \geq s,\\
u(s)&=&\psi\in L^2(D), \nonumber\\
u(t)|_{\partial D}&=&0\nonumber.
\end{eqnarray}
Here ${\cal L}$ is the second order differential operator with Dirichlet boundary condition on $D$,
\begin{eqnarray}\label{1.5}
{\cal L}u={1\over 2} \sum_{i,j=1}^d {\partial \over \partial x_j} \left(a_{ij}(x){\partial u \over {\partial x_i}}\right)+c(x)u.
\end{eqnarray}
Assume \\
{\bf Condition (L)}: the coefficients $a_{ij}, c$ are smooth functions on $\bar D$, $a_{ij}=a_{ji}$, and there exists $\gamma>0$ such that $\sum_ {i,j=1}^d a_{ij} \xi_i\xi_j\geq \gamma |\xi|^2$ for any $\xi=(\xi_1, \xi_2,\cdots,\xi_d)\in R^d$. \\
 \\
Under the above conditions, ${\cal L}$ is a self-adjoint uniformly elliptic operator and has discrete real-valued eigenvalues $\mu_1\geq \mu_2\geq \cdots$ such that $\mu_k\to{-\infty}$ when $k\to \infty$. Denote by $\{\phi_k\in L^2(D),\ k\geq 1\}$ a complete orthonormal system of  eigenfunctions of ${\cal L}$ with corresponding eigenvalues $\mu_k,\ k\geq 1$. Here the space $L^2(D)$ is a standard square integrable measurable function space vanishing on the boundary with norm $||\cdot||_{L^2(D)}$. A standard notation $H_0^1(D)$ denotes a standard Sobolev space of the square integrable measurable functions having the first order weak derivative in $L^2(D)$ and vanishing at the boundary $\partial D$. This is a Hilbert space with inner product $(u,v)=\int_D u(x)v(x)dx+\int_D (Du(x), Dv(x))dx$, for any $u, v\in H_0^1(D)$.
   From the uniformly elliptic condition, it's not difficult to know that $\phi_k\in H_0^1(D)$ and there exists a constant $C$ such that 
\begin{eqnarray}\label{eqn1.6}
||\nabla \phi_k||_{L^2(D)}\leq C\sqrt {|\mu_k|}.
\end{eqnarray}
 We will use it in the proof of our main theorem.
 
We assume the driving noise 
$W^k$ are mutually independent one-dimensional two-sided standard Brownian motions on the probability space ($\Omega, \cal F, P$) and $\sum_{k=1}^\infty\sigma_k^2(t)<\infty$.
Denote $\Delta:=\{(t,s)\in R^2, s\leq t\}$. Equation (\ref {march17a}) generates
a semi-flow $u:\Delta\times H\times\Omega\to H$ when the
solution exists uniquely in the space $H=L^2(D)$. Define $\theta: (-\infty,\infty)\times\Omega\to \Omega$ by $\theta_t\omega^k(s)=W^k(t+s)-W^k(t)$. Therefore ($\Omega, {\cal F}, P, (\theta_t)_{t\in R}$) is a metric dynamical system. Function $F:R\times R \to R$ is a continuous function. Without causing confusion of notation, we define Nemytskii operator $F: R\times L^2(D)\to L^2(D)$ with the same notation
$$F(t, u(t))(x)=F(t, u(t,x)),\ F^i(t, u(t))(x)=\int_D F(t, u(t))(y) \phi_i(y)dy\phi_i(x),\ x\in D,\ u\in L^2(D).$$
Assume $F$ and $\sigma_k$ satisfy: \vskip5pt

\noindent {\bf Condition (P)} {\it There exists a constant $\tau>0$
such that for any $t\in R$, $u\in L^2(D)$
\begin{eqnarray*}
F(t,u)=F(t+\tau,u),\
\sigma_k(t)=\sigma_k(t+\tau).
\end{eqnarray*}}

First,  we give the definition of the random periodic solution
\begin{defi}
A random periodic solution of period $\tau$ of a semi-flow $u:
\Delta\times L^2(D)\times\Omega\to L^2(D)$ is an ${\cal F}$- measurable
map $\varphi:(-\infty, \infty)\times \Omega\to L^2(D)$ such that
\begin{eqnarray}
u(t+\tau, t, \varphi(t,\omega), \omega)=\varphi(t+\tau,\omega)=\varphi(t, \theta_\tau \omega),
\end{eqnarray}
for any $t\in R$ and $\omega\in \Omega$.
\end{defi}

Instead of following the traditional geometric method of establishing the Poincar\'e mapping and finding its fixed point, in this paper,  we will push the new analysis method of coupled infinite horizon forward-backward
integral equations to the stochastic partial differential equations. This is the first paper dealing with the important question of periodic solution to stochastic partial differential equations.

We apply our result to the perturbation problem (\ref{1.1}) and (\ref{1.2}) we posed in the case when $H=R^d$, and the case when $H=L^2(D)$, $A={\cal L}$ a second order differential operator (\ref{1.5}) on a smooth bounded domain $D$. Assume the deterministic system has a periodic solution $Z$ which is hyperbolic. Denote by $G$ the graph of the periodic solution in $H$. Let $N$ be large enough such that the open ball with center $0$ and radius $N$ covers $G$. One can then define a differentiable function (assuming $f$ is differentiable) such that
\begin{eqnarray*}
f_N(u)=\chi({||u||^2\over N^2})f(u).
\end{eqnarray*}
Here  $\chi: R^1\to R^1$ is a smooth function such that 
\begin{eqnarray*}
\chi (z)=\left\{
\begin{array}{cl}
1, &{\rm when} \ |z|\leq 1,\\
0, &  {\rm when} \ |z|\ge 4.
\end{array} 
\right.
\end{eqnarray*}
It is easy to see that the truncated system 
\begin{eqnarray}\label{truncated}
{{du}\over {dt}} =Au+f_N(u)
\end{eqnarray}
has the same periodic solution $Z$ as Equation (\ref{1.1}). Our results imply that the perturbed system to Equation (\ref{truncated}) by an additive noise considered in \cite{feng-zhao-zhou} and in this paper respectively has a random periodic solution.

\section{Forward-backward infinite horizon stochastic integral equations}\label{section3}
\setcounter{equation}{0}



We consider the semilinear stochastic partial
differential equation (\ref{march17a}).
Denote the solution by $u(t,s,\omega, x)$.
Throughout this paper, we suppose that $\cal L$ is hyperbolic, i.e. none of the eigenvalues of ${\cal L}$ is zero, and $T_t=e^{{\cal L}t}$ is a hyperbolic linear flow induced by ${\cal L}$. 
 So $L^2(D)$ has a direct sum decomposition:
$$L^2(D)=E^s\oplus E^u,$$
where $$E^s=span\{v: v {\rm \  is\ a \ generalized\  eigenvector\  for\  an\  eigenvalue} \  \mu {\rm \ with \ }\mu<0\},$$
$$E^u=span\{v: v {\rm \  is\ a \ generalized\  eigenvector\  for\  an\  eigenvalue} \  \mu {\rm \ with \ }\mu>0\}.$$
Denote $\mu_m$ is the smallest positive eigenvalue of $\cal L$, and $\mu_{m+1}$ is the largest negative one.
We also define the projections onto
each subspace by
$$P^+:L^2(D)\rightarrow E^u, \ P^-: L^2(D) \rightarrow E^s.$$
Define ${\cal F}_s^t:=\sigma(W_u-W_v, s\leq v\leq u\leq t)$ and ${\cal F}^t:=\vee_{s\leq t}{\cal F}_s^t$.  
  The solution of the initial value problem (\ref{march17a}) is given by the following variation of constant formula:
  \begin{eqnarray}
  u(t,s,\psi,\omega)(x)&=&T_{t-s}\psi(x)+\int _s^t T_{t-r}F(r,u(r,s,\psi,\omega))(x)dr+\sum_{k=1}^\infty \int _s^t\sigma_k(r)(T_{t-r}\phi_k)(x)dW^k(r)\nonumber\\
  &=&\int_DK(t-s,x,y)\psi(y)dy+\int_s^t\int_DK(t-r,x,y)F(s,u(r,s,\psi,\omega))(y)dydr\nonumber\\
  &&+\sum_{k=1}^\infty \int_s^t\int_D K(t-r,x,y)\sigma_k(r)\phi_k(y)dydW^k(r), 
  \end{eqnarray}
where $K(t,x,y)$ is the heat kernel of the second order differential operator $\cal L$, 
\begin{eqnarray*}
(T_t\phi)(x)=\int_DK(t,x,y)\phi(y)dy,
\end{eqnarray*}
defines a linear operator $T_t:L^2(D)\to L^2(D)$ and $\int _s^t\sigma_k(r)(T_{t-r}\phi_k)(\cdot)dW^k(r)$ is an $L^2(D)$-valued stochastic integral.
Because ${\cal L}$ is a compact self-adjoint operator under the condition of this paper, so by Mercer's theorem (Chapter 3, Theorem 17, \cite{hoch}), we have
\begin{eqnarray*}
K(t, x,y)=\sum_{i=1}^\infty e^{\mu_i t}\phi_i(x)\phi_i(y).
\end{eqnarray*}
We consider a solution of the following coupled forward-backward infinite horizon stochastic integral
equation, which is a ${\cal B} (R)\otimes{\cal B} (D)\otimes\mathcal{F}$-measurable map
$Y: (-\infty,\infty)\times\Omega\rightarrow L^2(D)$ satisfying
\begin{eqnarray}\label{sep17a}
\hskip-0.7cmY(t,\omega)&=&\int _{-\infty}^tT_{t-s}P^-F(s,Y(s,\omega))ds-\int
_t^{\infty} T_{t-s}P^+F(s,Y(s,\omega))ds\nonumber
\\
&&+ (\omega)\sum_{k=1}^\infty \int _{-\infty}^t\sigma_k(s)T_{t-s}P^-\phi_k\,dW^k(s)-
(\omega)\sum_{k=1}^\infty \int _t^{\infty}\sigma_k(s)T_{t-s}P^+\phi_k\,dW^k(s)
\end{eqnarray}
for all $\omega\in \Omega$, $t\in(-\infty,\infty)$. The value of $Y(t,\omega)\in L^2(D)$ at $x$ is $Y(t,\omega)(x)$. Sometimes we write as $Y(t,\omega,x)$ when there is no confusing. We will give the following general theorem
which identifies the solution of the equation (\ref{sep17a}) and a random periodic
solution of stochastic differential equation (\ref{march17a}). First, we recall the
definition of a tempered random variable (Definition 4.1.1 in \cite{ar}):
\begin{defi}
A random variable $X: \Omega \to L^2(D)$ is called tempered with respect to the dynamical system $\theta$ if
$$
\lim\limits _{r\to \pm \infty}{1\over |r|}\log ||X(\theta _r\omega)||_{L^2(D)}= 0.
$$
The random variable is called tempered from above (below) if in the above limit,
the function $\log$ is replaced by $\log^+$ ($\log^-$), the positive (negative)
part of the function $\log$.
\end{defi}

\begin{thm}\label{aug20d} Assume Condition (P). If Cauchy problem (\ref{march17a}) has a unique
solution $u(t,s,\omega, x)$ and the coupled forward-backward infinite horizon stochastic integral equation (\ref{sep17a}) has
one solution $Y: (-\infty,+\infty)\times \Omega\rightarrow L^2(D)$ such that $Y(t+\tau,\omega)=Y(t,\theta_{\tau} \omega) \ {\rm for \ any} \
t\in R$ a.s., then $Y$ is a random periodic solution of
equation (\ref{march17a}) i.e.
\begin{eqnarray}
u(t+\tau,t, Y(t,\omega),\omega)=Y(t+\tau,\omega)=Y(t,\theta_{\tau} \omega) \ \ {\rm for \ any} \ \
t\in R
\ \ \ \ a.s.
\end{eqnarray}
Conversely, if equation (\ref{march17a}) has a random periodic solution $Y: (-\infty,+\infty)\times\Omega\rightarrow L^2(D)$
of period $\tau$ which is tempered from above for each $t$,
then $Y$ is a solution of the coupled forward-backward infinite horizon
stochastic integral equation (\ref{sep17a}).
\end{thm}
{\bf Proof:} Similar to the proof of Theorem 2.1 in \cite{feng-zhao-zhou}.
\hfill \hfill $\sharp$
\\

We will need the following generalized Schauder's fixed point theorem to prove our theorem. The proof was refined from the proof of Schauder's fixed point theorem and was given in \cite{feng-zhao-zhou}. 
\begin{thm} (Generalized Schauder's fixed point theorem)\label{Schauder}
Let $H$ be a Banach space, S be a convex subset of $H$. Assume a map $T: H\to H$ is continuous and $T(S)\subset S$
is relatively compact in $H$. Then $T$ has a fixed point in $H$.
\end{thm}

The generalized Schauder's fixed point theorem requires us to check the relative compactness. Since the equation can be transformed to an $\omega$-wise equation, one could be tempted to treat $\omega$ as a parameter and to try to define $\omega$-parameterised Banach space and subspace, and then to use Rellich-Kondrachov compactness embedding theorem to check the relative compactness. The problem with this approach is that, we get one solution with a parameter $\omega_1$ and one solution with a parameter $\omega_2$, but no priori relation between these solutions may be known. They may indeed belong to two different families of random periodic solutions due to the non-uniqueness of the solutions of the infinite horizon integral equation. Assume $\omega_2=\theta_\tau \omega_1$. It is desirable to have $Y(t+\tau, \omega_1)=Y(t, \omega_2)$ for all $t\geq 0$. But this is beyond what the analytic method can  offer to us immediately. To overcome this difficulty, we use Malliavin calculus, Wiener-Sobolev compact embedding theorem to get the relatively compactness of a sequence in $C^0([0, T], L^2(\Omega\times D))$ with Sobolev norm being bounded in $L^2(\Omega)$ and Malliavin derivative being bounded and equicontinuous in $L^2(\Omega\times D)$ uniformly in time. 
We denote by $C_p^\infty(R^n)$ the set of infinitely differentiable functions $f:R^n\to R$ such that $f$ and all its partial derivatives have polynomial growth. Let ${\mathcal S}$ be the class of smooth random variables $F$ such that $F=f(W(h_1),\cdots, W(h_n))$ with $n\in N$,  $h_1,\cdots, h_n\in L^2([0,T])$ and $f\in C_p^\infty(R^n)$, $W(h_i)=\int_0^T h_i(s)dW(s)$. The derivative operator of a smooth random variable $F$ is the stochastic process $\{{\cal D}_t F,\  t\in [0,T]\}$ defined by (c.f. \cite{nuallart})
$${\cal D}_t F=\sum_{i=1}^n {{\partial f}\over {\partial x_i}}(W(h_1),\cdots,W(h_n))h_i(t).$$
We will denote ${\cal D}^{1,2}$ the domain of ${\cal D}$ in $L^2(\Omega)$, i.e. ${\cal D}^{1,2}$ is the closure of ${\mathcal S}$ with respect to the norm
 $$||F||_{1,2}^2=E|F|^2+E||{\cal D}_tF||^2_{L^2([0,T])}.$$
 Denote $C^0([0,T],L^2(\Omega\times D))$ the set of continuous functions $f(\cdot,\cdot, \omega)$ with the norm
 $$||f||^2=\sup_{t\in [0,T]} \int_DE|f(t,x)|^2dx<\infty.$$
 It's easy to check the following refined version of relative compactness of Wiener-Sobolev space in Bally-Saussereau \cite{bally} also holds. This kind of compactness as a purely random variable
 version without including time and space variables was investigated by Da Prato, Malliavin and Nualart \cite{da-mall}  and Peszat \cite{peszat} first. Bally-Saussereau considered the convergence in $L^2([0,T]\times \Omega\times D)$. But the convergence in $L^2([0,T]\times \Omega\times D)$ is not enough for us in this paper. We consider the convergence in $C^0([0,T],L^2(\Omega\times D))$. Feng, Zhao and Zhou \cite{feng-zhao-zhou} used the compactness of a sequence of stochastic processes in $C^0([0,T], L^2(\Omega))$ to study periodic solution of stochastic differential equations. The Wiener-Sobolev compact embedding provides a powerful method to study the convergence of a sequence of random fields. This is a new direction of Malliavin calculus. The traditional application of Malliavin calculus was in regularity of densities and was studied intensively in literature.
  \begin{thm}\label{B-S}
Let $D$ be a bounded domain in $R^d$. Consider a sequence $(v_n)_{n\in N}$ of $C^0([0,T],L^2(\Omega\times D))$. Suppose that:\\
(1) $\sup_{n\in N}\sup_{t\in [0, T]} E||v_n(t, \cdot)||_{H^1(D)}^2 <\infty$.\\
(2) $\sup_{n\in N}\sup_{t\in [0, T]} \int_D ||v_n(t, x,\cdot)||_{1,2}^2 dx<\infty$.\\
(3) There exists a constant $C>0$ such that for any $t_1, t_2\in [0, T]$ 

         $\hskip 0.5cm \sup_n\int_D E|v_n(t_1,x)-v_n(t_2,x)|^2dx< C|t_1-t_2|.$
%
\\
 (4) (4i) There exists a constant $C$ such that for any $0<\alpha<\beta<T$, and $h \in R$ with $|h|<\min (\alpha, T-\beta)$,
 
 $\hskip 0.3cm$  and any $t_1, t_2\in [0,T]$, 
 
  $\hskip 0.5cm \sup_n\int_D\int_{\alpha}^{\beta}E |{\cal D}_{\theta+h}v_n(t_1,x)-{\cal D}_\theta v_n(t_2,x)|^2 d\theta dx<C(|h|+|t_1-t_2|)$.
 
 (4ii) For any $\epsilon>0$, there exist $0<\alpha<\beta<T$ such that
 
$\hskip 0.5cm\sup_n \sup_{t\in[0, T]}\int_D\int_{[0, T]\backslash (\alpha,\beta)} E|{\cal D}_\theta v_n(t,x)|^2 d\theta dx <\epsilon$.\\
Then $\{v_n, n\in N\}$ is relatively compact in $C^0([0,T],L^2(\Omega\times D))$.
 \end{thm}
{\bf Proof:} Recall the Wiener chaos expansion 
$$v_n(t,\omega, x)=\sum_{m=0}^\infty I_m(f_n^m(\cdot,t,x))(\omega),$$
where $f_n^m(\cdot,t,x)$  are symmetric elements of $L^2([0, T]^m\times D)$ for each $m\geq 0$. When $m=0$, $f_n^0(t,x)=Ev_n(t,x)$, and
\begin{eqnarray*}
\sup_n||f_n^0(t,\cdot)||_{H^1(D)}^2\leq \sup_n E||v_n(t,\cdot)||_{H^1(D)}^2<\infty.
\end{eqnarray*}
So $f_n^0(t,x)$ is relatively compact in $L^2(D)$ for fixed $t\in [0,T]$ by Rellich-Kondrachov compact embedding theorem.
But for any $t_1, t_2\in [0, T]$,
\begin{eqnarray*}
&&\sup_n\sup_{t\in [0,T]} ||f_n^0(t,\cdot)||_{L^2(D)}^2\leq \sup_n\sup_{t\in [0,T]}{E||v_n(t,\cdot)||_{L^2(D)}^2}<\infty,\\
&&\sup_n||f_n^0(t_1,\cdot)-f_n^0(t_2,\cdot)||_{L^2(D)}^2\leq \sup_n {E||v_n(t_1)-v_n(t_2)||_{L^2(D)}^2}\leq {C|t_1-t_2|}.
\end{eqnarray*}
So by Arzela-Ascoli lemma, $\{f_n^0\}_{n=1}^\infty$ is relatively compact in $C^0([0, T],L^2(D))$.
  For each $m\geq 1$, using the same argument as in Bally-Saussereau \cite{bally}, we conclude for each fixed t, $\{f_n^m(\cdot, t,x)\}_{n\in N}$ is relatively compact in $L^2([0, T]^m\times D)$. Moreover, for each $t_1, t_2\in [0,T]$, consider
\begin{eqnarray*}
&&\sup_n\int_D||f_n^m(\cdot, t_1,x)-f_n^m(\cdot,t_2,x)||_{L^2([0, T]^m)}^2dx\\
&\leq& \sup_n\int_D\int_0^T E|{\cal D}_\theta v_n(t_1,x)-{\cal D}_\theta v_n(t_2,x)|^2d\theta dx\\
&\leq& C|t_1-t_2|,
\end{eqnarray*}
and
\begin{eqnarray*}
\sup_n\sup_{t\in [0, T]} \int_D||f_n^m(\cdot,t,x)||_{L^2([0, T]^m)}^2dx\leq\sup_n\sup_{t\in [0, T]}\int_D \int_0^TE|{\cal D}_\theta v_n(t,x)|^2d\theta dx<\infty.
\end{eqnarray*}
Then by Arzela-Ascoli lemma, we know that $\{f_n^m\}_{n=1}^\infty$ is relatively compact in $C^0([0, T], L^2([0, T]^m\times D))$. Thus we can conclude $\{v_n\}_{n=1}^\infty$ is relatively compact in $C^0([0, T], L^2(\Omega\times D))$ using the same argument as in \cite{bally}.
\hfill \hfill $\sharp$\\

Now we are going to prove that equation
(\ref{sep17a}) has a solution under some conditions. So according to Theorem \ref{aug20d},
this gives the existence of the random periodic solution for the stochastic evolution equation (\ref{march17a}).

\begin{thm}\label{aug20b}
Assume the coefficients of the second order differential operator ${\cal L}$ 
satisfy condition (L) and the operator ${\cal L}$ is hyperbolic. Let $F:(-\infty, \infty) \times R\to R$ be a continuous map, globally bounded and $\nabla F(t,\cdot)$ being globally bounded,
    and $F$ and $\sigma_k$ also satisfy Condition (P) and $\sum_{k=1}^\infty|\sigma_k(t)|^2<\infty$, and there exists a constant $L_1>0$ such that $\sum_{k=1}^\infty|\sigma_k(s_1)-\sigma_k(s_2)|^2\leq L_1 |s_1-s_2|$. Then there exists at least one ${\cal B}(R)\otimes\mathcal{F}$-measurable map
$Y: (-\infty,+\infty)\times\Omega\rightarrow L^2(D)$ satisfying equation (\ref{sep17a})
 and  $Y(t+\tau,\omega)=Y(t, \theta_\tau\omega)$ for
any $t\in R$, $\omega\in \Omega$.
\end{thm}

The proof of the theorem is very complex and is based on the following observation and a series of lemmas.
Define the ${\cal B}(R)\otimes\mathcal{F}$-measurable map
$Y_{1}:(-\infty,+\infty)\times\Omega\rightarrow L^2(D)$ by
\begin{eqnarray}
\hskip -0.7cmY_{1}(t,\omega)=(\omega)\sum_{k=1}^\infty\int _{-\infty}^t\sigma_k(s)T_{t-s}P^-\phi_k\,dW^k(s)-
(\omega)\sum_{k=1}^\infty\int _t^{\infty}\sigma_k(s)T_{t-s}P^+\phi_k\,dW^k(s).
\end{eqnarray}
Then by changing of variable and periodicity of $\sigma_k$, we have
\begin{eqnarray}
&&Y_{1}(t,\theta_\tau\omega)\nonumber\\
&=&(\theta_\tau\omega)\sum_{k=1}^\infty\int^{t}_{-\infty}\sigma_k(s)T_{t-s}P^-\phi_k\,dW^k(s)-(\theta_\tau\omega)\sum_{k=1}^\infty\int^{\infty}_{t}\sigma_k(s)T_{t-s}P^+\phi_k\,dW^k(s)\nonumber\\
&=&(\omega)\sum_{k=1}^\infty\int^{t+\tau}_{-\infty}\sigma_k(s)T_{t+\tau-s}P^-\phi_k\,dW^k(s)-(\omega)\sum_{k=1}^\infty\int^{\infty}_{t+\tau}\sigma_k(s)T_{t+\tau-s}P^+\phi_k\,dW^k(s)\nonumber\\
&=&Y_1(t+\tau,\omega).
\end{eqnarray}
On the other hand, 
\begin{eqnarray*}
Y_{1}(t,\omega, x)
&=& \sum_{k=1}^\infty\sum_{i=m+1}^\infty\int^{t}_{-\infty}  e^{\mu_i (t-s)}\sigma_k(s)\int_D\phi_i(y)\phi_k(y)\,dy\phi_i(x)dW^k(s)\\
&&-\sum_{k=1}^\infty\sum_{i=1}^m\int_{t}^{\infty}  e^{\mu_i (t-s)}\sigma_k(s)\int_D\phi_i(y)\phi_k(y)\,dy\phi_i(x)dW^k(s)\\
&=& \sum_{i=m+1}^\infty \int_{-\infty}^t e^{\mu_i (t-s)}\sigma_i(s)dW^i(s)\phi_i(x)-\sum_{i=1}^m \int_t^\infty e^{\mu_i (t-s)}\sigma_i(s)dW^i(s)\phi_i(x),
\end{eqnarray*}
as $\{\phi_i\}$ is the basis of $L^2(D)$, so $\int_D\phi_i(y)\phi_j (y) dy=0$, when $i\neq j$ and $\int_D\phi_i^2(y)=1$.
Moreover, we can calculate
\begin{eqnarray*}
||Y_1||^2
&=&\sup_t E\int_D| Y_1(t,y)|^2dy\\
&\leq & 2 \sup_t E\int_D\Big|\sum_{i=m+1}^\infty \int_{-\infty}^t e^{\mu_i(t-s)} \sigma_i(s)dW^i(s)\phi_i(y)\Big|^2dy\\
&&+2 \sup_t E\int_D\Big|\sum_{i=1}^m \int^{\infty}_t e^{\mu_i(t-s)}\sigma_i(s)dW^i(s)\phi_i(y)\Big|^2dy\\
&=& 2\sup_t E\int_D\sum_{i=m+1}^\infty \int_{-\infty}^t e^{2\mu_i(t-s)}|\sigma_i(s)|^2ds| \phi_i(y)|^2dy\\
&&+2 \sup_t E\int_D\sum_{i=1}^m \int^{\infty}_t e^{2\mu_i(t-s)}|\sigma_i(s)|^2ds| \phi_i(y)|^2dy\\
&\leq& 2\sup_t E\sum_{i=m+1}^\infty \int_{-\infty}^t e^{2\mu_{m+1}(t-s)}|\sigma_i(s)|^2ds\\
&&+2\sup_t E\sum_{i=1}^m \int^{\infty}_t e^{2\mu_m(t-s)}|\sigma_i(s)|^2ds\\
&\leq& (-{1\over{\mu_{m+1}}}+{1\over {\mu_m}})\sup_{s\in (-\infty,\infty)}\sum_{i=1}^\infty \sigma_i^2(s)\\
&<&\infty.
\end{eqnarray*}

Secondly, we need to solve the equation
\begin{eqnarray}\label{zhao5}
Z(t,\omega)&=&\int^{t}_{-\infty}T_{t-s}P^-F(s, Z(s,\omega)+Y_1(s,\omega)))ds\nonumber\\
&&-\int^{\infty}_{t}T_{t-s}P^+F(s,Z(s,\omega)+Y_1(s,\omega)))ds.
\end{eqnarray}
For this we define
 \begin{eqnarray*}
 &&C_{\tau}^0((-\infty, +\infty), L^2(\Omega\times D))\nonumber\\
 &:=&\{f\in C^0((-\infty, +\infty), L^2(\Omega\times D)):
 \ {\rm for \ any} \ \
t\in (-\infty,\infty),\
f(\tau+t,\omega,x)=f(t,\theta_{\tau}\omega,x) \},
\end{eqnarray*}
and for any $z\in  C_{\tau}^0((-\infty, +\infty), L^2(\Omega\times D))$, define
\begin{eqnarray}\label{2.7}
&&{\cal M}(z)(t,\omega,x)\nonumber\\
&:=&\int^{t}_{-\infty}T_{t-s}P^-F(s, z(s,\omega)+Y_1(s,\omega))(x)ds-\int^{\infty}_{t}T_{t-s}P^+F(s,z(s,\omega)+Y_1(s,\omega))(x)ds.
\end{eqnarray}
The idea is to find a fixed point to ${\cal M}$ in $C_{\tau}^0((-\infty, +\infty), L^2(\Omega\times D))$ using the generalized Schauder's fixed point Theorem \ref {Schauder}.
\begin{lem}\label{lem2.1}
Under the conditions of Theorem \ref{aug20b}, the map 
$${\cal M}: C_{\tau}^0((-\infty, +\infty), L^2(\Omega\times D))\to C_{\tau}^0((-\infty, +\infty), L^2(\Omega\times D))$$ is a continuous map. Moreover ${\cal M}$ maps $C_{\tau}^0((-\infty, +\infty), L^2(\Omega\times D))$ into $C_{\tau}^0((-\infty, +\infty), L^2(\Omega\times D))\cap L^\infty((-\infty, +\infty), L^2(\Omega,H^1_0(D)))$. 
\end{lem}
{\bf Proof:} Firstly, for any $z\in C_{\tau}^0((-\infty, +\infty), L^2(\Omega\times D))$, from $\{\phi_i\}$ is the basis of $L^2(D)$, 
Cauchy-Schwarz inequality and the linear growth of $F$ with respect to the second variable, we have 
\begin{eqnarray*}
&&E\int_D|{\cal M}(z)(t,x)|^2dx\\
&\leq & 2\int_DE\Big|\int_{-\infty}^t \int_D\sum_{i=m+1}^\infty e^{\mu_i(t-s)}\phi_i(x)\phi_i(y)F^i(s, z(s)+Y_1(s))(y)dyds\Big|^2dx\\
&&+2\int_DE\Big|\int_t^{\infty} \int_D\sum_{i=1}^m e^{\mu_i(t-s)}\phi_i(x)\phi_i(y)F^i(s, z(s)+Y_1(s))(y)dyds\Big|^2dx\\
&=& 2E\sum_{i=m+1}^\infty\Big|\int_{-\infty}^t\int_De^{\mu_i(t-s)}\phi_i(y)F^i(s, z(s)+Y_1(s))(y)dyds\Big|^2\\
&&+2E\sum_{i=1}^m\Big|\int_t^\infty \int_De^{\mu_i(t-s)}\phi_i(y)F^i(s, z(s)+Y_1(s))(y)dyds\Big|^2\\
&\leq&2E\sum_{i=m+1}^\infty\Big[\int^t_{-\infty} \int_D e^{\mu_i(t-s)}|\phi_i(y)|^2dyds\cdot\int_{-\infty}^t \int_De^{\mu_i(t-s)}|F^i(s, z(s)+Y_1(s))(y)|^2dyds\Big]\\
&&+2E\sum_{i=1}^m\Big[\int_t^{\infty} \int_D e^{\mu_i(t-s)}|\phi_i(y)|^2dyds\cdot\int^{\infty}_t \int_De^{\mu_i(t-s)}|F^i(s, z(s)+Y_1(s))(y)|^2dyds\Big]\\
&\leq&(-{2\over{\mu_{m+1}}})\sum_{i=m+1}^\infty E\int_{-\infty}^t \int_De^{\mu_{m+1}(t-s)}|F^i(s, z(s)+Y_1(s))(y)|^2dyds\\
&&+{2\over{\mu_{m}}}\sum_{i=1}^mE\int_t^{\infty} \int_De^{\mu_{m}(t-s)}|F^i(s, z(s)+Y_1(s))(y)|^2dyds\\
&\leq&2||F||^2_\infty(\frac{1}{\mu^2_{m+1}}+\frac{1}{\mu^2_{m}})vol(D)\\
&<&\infty.
\end{eqnarray*}
We prove that ${\cal M} (z)(\cdot, \omega,x)$ is continuous in $L^2(\Omega\times D)$, for $z\in C_{\tau}^0((-\infty, +\infty), L^2(\Omega\times D))$.
 For this, taking any $t_1$, $t_2 \in
(-\infty, +\infty)$ with $t_1\leq t_2$, we have
\begin{eqnarray*}
&&E \int_D|{\cal M}(z)({t_1},x)-
{\cal M}(z)({t_2},x)|^2dx\\
&\leq& 2\int_DE\Big [\big|
\int^{t_1}_{-\infty}T_{t_1-s}P^-F(s,z({s})+Y_{1}({s}))(x)ds
-\int^{t_2}_{-\infty}T_{t_2-s}P^-F(s,z({s})+Y_{1}({s}))(x)ds\big|^2\\
&&+\big|
\int^{+\infty}_{t_2}T_{t_2-s}P^+F(s,z({s})+Y_{1}({s}))(x)ds-
\int^{+\infty}_{t_2}T_{t_2-s}P^+F(s,z({s})+Y_{1}({s}))(x)ds
\big|^2\Big]dx.
\end{eqnarray*}
For the first term, considering $\{\phi_i\}$ is the basis of $L^2(D)$, and noting the following simple computation, for $i\geq m+1$, 
\begin{eqnarray*}
\int_{-\infty}^{t_1}|e^{\mu_i(t_1-s)}-e^{\mu_i(t_2-s)}|ds=\int_{-\infty}^{t_1}e^{\mu_i(t_1-s)}|1-e^{\mu_i(t_2-t_1)}|ds
\leq(t_2-t_1)\int_{-\infty}^{t_1}e^{\mu_i(t_1-s)} |\mu_i| ds=t_2-t_1,
\end{eqnarray*}
we have the following estimate,
\begin{eqnarray*}
&&\int_DE\big|
\int^{t_1}_{-\infty}T_{t_1-s}P^-F(s,z({s})+Y_{1}({s}))(x)ds
-\int^{t_2}_{-\infty}T_{t_2-s}P^-F(s,z({s})+Y_{1}({s}))(x)ds
\big|^2dx \\
&\leq&2 \int_DE\big|\int_{-\infty}^{t_1}\int_D\sum_{i=m+1}^\infty(e^{\mu_i(t_1-s)}-e^{\mu_i(t_2-s)})\phi_i(x)\phi_i(y)F^i(s,z(s)+Y_1(s))(y)dyds\big|^2dx\\
&&+2 \int_DE\big|\int_{t_1}^{t_2}\int_D\sum_{i=m+1}^\infty e^{\mu_i(t_2-s)}\phi_i(x)\phi_i(y)F^i(s,z(s)+Y_1(s))(y)dyds\big|^2dx\\
&\leq& 2E\sum_{i=m+1}^\infty\big|\int_{-\infty}^{t_1}\int_D(e^{\mu_i(t_1-s)}-e^{\mu_i(t_2-s)})\phi_i(y)F^i(s,z(s)+Y_1(s))(y)dyds\big|^2\\
&&+2E\sum_{i=m+1}^\infty\big|\int_{t_1}^{t_2}\int_D e^{\mu_i(t_2-s)}\phi_i(y)F^i(s,z(s)+Y_1(s))(y)dyds\big|^2\\
&\leq& 2E\sum_{i=m+1}^\infty \big[\int_{-\infty}^{t_1}\int_D(e^{\mu_i(t_1-s)}-e^{\mu_i(t_2-s)})|\phi_i(y)|^2dyds\\
&&\hskip 3cm \cdot\int_{-\infty}^{t_1}\int_D(e^{\mu_i(t_1-s)}-e^{\mu_i(t_2-s)})|F^i(s,z(s)+Y_1(s))(y)|^2dyds\big]\\
&&+2E\sum_{i=m+1}^\infty\int_{t_1}^{t_2}\int_D|\phi_i(y)|^2dyds\cdot\int_{t_1}^{t_2}\int_D|F^i(s,z(s)+Y_1(s))(y)|^2dyds\\
&\leq& 2E\sum_{i=m+1}^\infty (t_2-t_1)\int_{-\infty}^{t_1}\int_De^{\mu_{m+1}(t_1-s)}| F^i(s,z(s)+Y_1(s))(y)|^2dyds\\
&&+2(t_2-t_1)^2 ||F||^2_\infty vol(D)\\
&\leq&(-{2\over{\mu_{m+1}}})|t_2-t_1| ||F||^2_\infty vol(D)+2(t_2-t_1)^2 ||F||^2_\infty vol(D)\\
&\leq&C|t_2-t_1|.
\end{eqnarray*}
And by a similar argument to the second part, we have
\begin{eqnarray*}
&&E\big| \int^{+\infty}_{t_1}T_{t_1-s}P^+F(s,z(s) +
Y_{1}({s}))ds-
\int^{+\infty}_{t_2}T_{t_2-s}P^+F(s,z(s)+Y_{1}({s}))ds\big|^2 
\leq C|t_2-t_1|.
\end{eqnarray*}
Therefore, by combining two parts, we have
\begin{eqnarray*}
E\int_D| {\cal M}(z)({t_2,x})-
{\cal M}(z)({t_1,x})|^2dx\leq C|t_2-t_1|.
\end{eqnarray*}
Therefore we have ${\cal M}$ also maps  $C_{\tau}^0((-\infty, +\infty), L^2(\Omega\times D))$ into itself.
To
see the continuity, for any $z_1,z_2\in
C_{\tau}^0((-\infty, +\infty), L^2(\Omega\times D))$,
\begin{eqnarray*}
&&\int_DE
|{\cal M}(z_1)(t,x)-{\cal M}(z_2)(t,x)|^2dx\\
&\leq&2E\sum_{i=m+1}^\infty\Big[\int_{-\infty}^t\int_De^{\mu_i(t-s)}\phi_i(y)(F^i(s,z_1(s)+Y_1(s))(y)\\
&&\hskip 3cm-F^i(s,z_2(s)+Y_1(s))(y))dyds\Big]^2\\
&&+2E\sum_{i=1}^m\Big[\int_{-\infty}^t\int_De^{\mu_i(t-s)}\phi_i(y)(F^i(s,z_1(s)+Y_1(s))(y)\\
&&\hskip 3cm-F^i(s,z_2(s)+Y_1(s))(y))dyds\Big]^2\\
&\leq& 2 E \sum_{i=m+1}^\infty\Big[\int_{-\infty}^t\int_De^{\mu_i(t-s)}|\phi_i(y)|^2dyds\\
&&\hskip 1cm\cdot\int_{-\infty}^t\int_De^{\mu_i(t-s)}|F^i(s,z_1(s)+Y_1(s))(y)-F^i(s,z_2(s)+Y_1(s))(y)|^2dyds\Big]\\
&&+ 2 E \sum_{i=1}^m\Big[\int^{\infty}_t\int_De^{\mu_i(t-s)}|\phi_i(y)|^2dyds\\
&&\hskip 1cm\cdot\int^{\infty}_t\int_De^{\mu_i(t-s)}|F^i(s,z_1(s)+Y_1(s))(y)-F^i(s,z_2(s)+Y_1(s))(y)|^2dyds\Big]\\
&\leq & 2(-{1\over {\mu_{m+1}}})\int_{-\infty}^t\int_De^{\mu_{m+1}(t-s)}|F(s,z_1(s)+Y_1(s))(y)-F(s,z_2(s)+Y_1(s))(y)|^2dyds\\
&&+2{1\over {\mu_{m}}}\int^{\infty}_t\int_De^{\mu_{m}(t-s)}|F(s,z_1(s)+Y_1(s))(y)-F(s,z_2(s)+Y_1(s))(y)|^2dyds\\
&\leq &2||\nabla F||^2_\infty(\frac{1}{\mu_{m+1}^2}+\frac{1}{\mu_{m}^2})\sup_{t\in (-\infty,+\infty)}\int_DE|z_1(t,x)-z_2(t,x)|^2dx,
\end{eqnarray*}
where
$$||\nabla F||^2_\infty:=\sup_{t\in(-\infty,\infty), u\in R}|\nabla F(t,u)|^2=\sup_{t\in(-\infty,\infty), u\in R}\sum_{i=1}^\infty|\nabla F^i(t,u)|^2.$$
That is to say that ${\cal M}: C_{\tau}^0((-\infty, +\infty), L^2(\Omega\times D))\to
C_{\tau}^0((-\infty, +\infty), L^2(\Omega\times D))$ is a continuous map.
Secondly, we need to prove ${\cal M}(z)\in L^\infty((-\infty,\infty), L^2(\Omega, H^1(D)))$ for $z\in C_{\tau}^0((-\infty, +\infty), L^2(\Omega\times D))$. Note
\begin{eqnarray*}
&&E\int_D|\nabla_x {\cal M}(z)(t,x)|^2dx\\
&\leq& 2E\int_D \Big|\int_{-\infty}^t\int_D\sum_{i=m+1}^\infty e^{\mu_i(t-s)}\nabla_x\phi_i(x)\phi_i(y)F^i(s,z(s)+Y_1(s))(y)dyds\Big|^2dx\\
&&+2E\int_D \Big|\int_{-\infty}^t\int_D\sum_{i=1}^m e^{\mu_i(t-s)}\nabla_x\phi_i(x)\phi_i(y)F^i(s,z(s)+Y_1(s))(y)dyds\Big|^2dx\\
&:=&A_1+A_2.
\end{eqnarray*}
For $A_1$, by Cauchy-Schwarz inequality and (\ref{eqn1.6}), we have
\begin{eqnarray*}
A_1&=&2E\int_D \Big[\sum_{i,j=m+1}^\infty \int_{-\infty}^t\int_De^{\mu_i(t-s)}\nabla_x\phi_i(x)\phi_i(y)F^i(s,z(s)+Y_1(s))(y)dyds\\
&&\hskip 2cm\cdot \int_{-\infty}^t\int_De^{\mu_j(t-s)}\nabla_x\phi_j(x)\phi_j(y)F^j(s,z(s)+Y_1(s))(y)dyds\Big]dx\\
&\leq & E \sum_{i,j=m+1}^\infty\left(\int_D |\nabla_x\phi_i(x)|^2dx\int_D|\nabla_x\phi_j(x)|^2dx\right)^{1\over 2}\\
&&\hskip 2cm\cdot \int_{-\infty}^t\int_De^{\mu_i(t-s)}|\phi_i(y)||F^i(s,z(s)+Y_1(s))(y)|dyds\\
&&\hskip 2cm\cdot \int_{-\infty}^t\int_De^{\mu_j(t-s)}|\phi_j(y)||F^j(s,z(s)+Y_1(s))(y)|dyds\\
&\leq & 2CE\Bigg[ \sum_{i=m+1}^\infty \left (\int_{-\infty}^t\int_De^{\mu_i(t-s)}|\mu_i|^{1\over 2}|\phi_i(y)||F^i(s,z(s)+Y_1(s))(y)|dyds\right)^2\\
&&\hskip1cm\cdot\sum_{j=m+1}^\infty \left (\int_{-\infty}^t\int_De^{\mu_j(t-s)}|\mu_j|^{1\over 2}|\phi_j(y)||F^j(s,z(s)+Y_1(s))(y)|dyds\right)^2\Bigg]^{1\over 2}\\
&\leq &2CE\Bigg[\sum_{i=m+1}^\infty\left(\int^t_{-\infty} \int_D e^{\mu_i(t-s)}|\mu_i||\phi_i(y)|^2dyds\cdot\int_{-\infty}^t \int_De^{\mu_i(t-s)}|F^i(s, z(s)+Y_1(s))(y)|^2dyds\right)\\
&&\cdot\sum_{j=m+1}^\infty\left(\int^t_{-\infty} \int_D e^{\mu_j(t-s)}|\mu_j||\phi_j(y)|^2dyds\cdot\int_{-\infty}^t \int_De^{\mu_j(t-s)}|F^j(s, z(s)+Y_1(s))(y)|^2dyds\right)\Bigg]^{1\over 2}\\
&\leq&2C\Bigg[\sum_{i=m+1}^\infty\int_t^{\infty} \int_De^{\mu_{m+1}(t-s)}|F^i(s, z(s)+Y_1(s))(y)|^2dyds\\
&&\hskip1cm\cdot\sum_{j=m+1}^\infty\int_t^{\infty} \int_De^{\mu_{m+1}(t-s)}|F^i(s, z(s)+Y_1(s))(y)|^2dyds\Bigg]^{1\over 2}\\
&\leq &2C||F||_\infty^2 ({-{1\over {\mu_{m+1}}}}) vol(D)\\
&<&\infty.
\end{eqnarray*}
Similarly, 
\begin{eqnarray*}
A_2\leq 2C ||F||_\infty^2 ({{1\over {\mu_{m}}}}) vol(D)<\infty.
\end{eqnarray*}
Therefore, we can see ${\cal M}$ maps $C^0_{\tau}((-\infty, +\infty), L^2(\Omega\times D))$ into $L^\infty((-\infty, +\infty), L^2(\Omega, H^1_0(D)))$.
$\hfill \hfill \sharp$
\vskip0.5cm
Now let us define a subset of $C_{\tau}^0((-\infty,+\infty), L^2(\Omega\times D))$ as follows:
\begin{eqnarray*}
&&C^0_{\tau,\alpha}((-\infty,+\infty),L^2(D, {\cal D}^{1,2})) \\
&:=&\{ f\in C_{\tau}^0((-\infty,+\infty), L^2(\Omega\times D)):\ f|_{[0,\tau)}\in C^0([0,\tau),L^2(D, {\cal D}^{1,2})),\\
&& i.e.\  ||f||^2=\sup_{t\in [0,\tau)}\int_D||f(t,x)||^2_{1,2}dx<\infty, {\rm \ and\  for \ any}\  t,r\in [0,\tau),\ i=0, \pm 1,\pm 2,\cdots \nonumber\\
&&\int_DE|{\cal D}_r f(t,\theta_{i\tau}\cdot, x)|^2dx\leq \alpha_r(t), \sup_{s,r_1, r_2 \in [0,\tau)}{{\int_DE|{\cal D}_{r_1} f(s,\theta_{i\tau}\cdot, x)-{\cal D}_{r_2} f(s,\theta_{i\tau}\cdot, x)|^2dx}\over {|r_1-r_2|}}<\infty\}.
\end{eqnarray*}
Here $\alpha_r(t)$ is the solution of integral equation (see page 324 in \cite{polyanin}) 
\begin{eqnarray}\label{eqn2.15}
\alpha_r(t)=A\int_{r-2\tau}^{r+2\tau} e^{-\beta |t-s|}\alpha_r(s)ds +B,
\end{eqnarray}
where
\begin{eqnarray*}
&&A=C||\nabla F||^2_\infty(-{1\over {\mu_{m+1}}}\sum_{i=0}^\infty e^{\mu_{m+1}i\tau}+{1\over {\mu_{m}}}\sum_{i=0}^\infty e^{-\mu_{m}i\tau}), \\
&&B=C||\nabla F||^2_\infty \sup_{s\in (-\infty,\infty)} \sum_{j=1}^\infty \sigma_j^2(s)({1\over {\mu^2_{m+1}}}+{1\over {\mu^2_{m}}}),\  \beta=min\{-\mu_{m+1},\mu_m\}. 
\end{eqnarray*}
This is a convex set.\\
\begin{lem}\label{lem2.2}
Under the conditions of Theorem \ref{aug20b}, ${\cal M}$ maps $C^0_{\tau,\alpha}((-\infty,+\infty),L^2(D, {\cal D}^{1,2}))$ into itself.
\end{lem}
{\bf Proof:} The Malliavin derivatives of $Y_1(t,\omega,x)$ and ${\cal M} (z)(t,\omega,x)$ can be calculated as:
\begin{eqnarray}
&&\hskip -1cm {\cal D}_rY_1(t, \omega,x)=\left\{
\begin{array}{ll}
\sum\limits_{i=m+1}^\infty e^{\mu_i(t-r)} \phi_i(x)\sigma_i(r),
 &{\rm if} \  r\leq t,\\
-\sum\limits_{i=1}^m e^{\mu_i(t-r)} \phi_i(x)\sigma_i(r), \ &{\rm if} \ r> t.
\end{array}
\right.\label{2.16}
\end{eqnarray}
When $r\leq t$, it is easy to see that 
\begin{eqnarray}
&& {\cal D}_r{\cal M}(z)(t, \omega,x)\\
&=&\sum_{i=m+1}^\infty\left(\int_{-\infty}^t\int_D  e^{\mu_i(t-s)}\phi_i(y)\nabla F^i(s, z(s,\omega)+Y_1(s,\omega))(y){\cal D}_rz(s,\omega,y)dyds\right)\phi_i(x)\nonumber\\
&&- \sum_{i=1}^m\left(\int_t^\infty\int_De^{\mu_i(t-s)}\phi_i(y)\nabla F^i(s, z(s,\omega)+Y_1(s,\omega))(y){\cal D}_rz(s,\omega,y)dyds\right)\phi_i(x)\nonumber\\
&&+\sum_{i=m+1}^\infty\left(\int_{-\infty}^r \int_D e^{\mu_i(t-s)}\phi_i(y)\nabla F^i(s, z(s,\omega)+Y_1(s,\omega))(y)\sum_{j=1}^m\Big(-e^{\mu_j(s-r)}\phi_j(y)\sigma_j(r)\Big)dyds\right)\phi_i(x)\nonumber\\
&&+\sum_{i=m+1}^\infty\left(\int_r^t \int_D e^{\mu_i(t-s)}\phi_i(y)\nabla F^i(s, z(s,\omega)+Y_1(s,\omega))(y)\sum_{j=m+1}^\infty\Big(e^{\mu_j(s-r)}\phi_j(y)\sigma_j(r)\Big)dyds\right)\phi_i(x)\nonumber\\
&&-\sum_{i=1}^m\left(\int_t^\infty \int_D e^{\mu_i(t-s)}\phi_i(y)\nabla F^i(s, z(s,\omega)+Y_1(s,\omega))(y)\sum_{j=m+1}^\infty\Big(e^{\mu_j(s-r)}\phi_j(y)\sigma_j(r)\Big)dyds\right)\phi_i(x).\nonumber
\end{eqnarray}
Similarly, when $r> t$, we have
\begin{eqnarray}
&& {\cal D}_r{\cal M}(z)(t, \omega)\\
&=&\sum_{i=m+1}^\infty\left(\int_{-\infty}^t\int_D e^{\mu_i(t-s)}\phi_i(y)\nabla F^i(s, z(s,\omega)+Y_1(s,\omega))(y){\cal D}_rz^j(s,\omega,y)dyds\right)\phi_i(x)\nonumber\\
&&-\sum_{i=1}^m\left(\int_t^\infty\int_D e^{\mu_i(t-s)}\phi_i(y)\nabla F^i(s, z(s,\omega)+Y_1(s,\omega))(y){\cal D}_rz^j(s,\omega,y)dyds\right)\phi_i(x)\nonumber\\
&&+\sum_{i=m+1}^\infty\left(\int_{-\infty}^t \int_D e^{\mu_i(t-s)}\phi_i(y)\nabla F^i(s, z(s,\omega)+Y_1(s,\omega))(y)\sum_{j=1}^m\Big(-e^{\mu_j(s-r)}\phi_j(y)\sigma_j(r)\Big)dyds\right)\phi_i(x)\nonumber\\
&&-\sum_{i=1}^m\left(\int_t^r \int_D e^{\mu_i(t-s)}\phi_i(y)\nabla F^i(s, z(s,\omega)+Y_1(s,\omega))(y)\sum_{j=1}^m\Big(-e^{\mu_j(s-r)}\phi_j(y)\sigma_j(r)\Big)dyds\right)\phi_i(x)\nonumber\\
&&-\sum_{i=1}^m\left(\int_r^\infty \int_D e^{\mu_i(t-s)}\phi_i(y)\nabla F^i(s, z(s,\omega)+Y_1(s,\omega))(y)\sum_{j=m+1}^\infty\Big(e^{\mu_j(s-r)}\phi_j(y)\sigma_j(r)\Big)dyds\right)\phi_i(x).\nonumber
\end{eqnarray}
So using Cauchy-Schwarz inequality, we have for any $k=0, \pm 1, \pm 2, \cdots$, $z\in C^0_{\tau,\alpha}((-\infty,+\infty),L^2(D, {\cal D}^{1,2}))$, 
when $0\leq r\leq t< \tau$,
\begin{eqnarray*}
&&E\int_D|{\cal D}_r{\cal M}(z)(t,\theta_{k\tau}\cdot, x)|^2dx\\
&\leq &CE\sum_{i=m+1}^\infty\Big[\int_{-\infty}^t\int_D e^{\mu_i(t-s)}\phi_i(y)\nabla F^i(s, z(s,\theta_{k\tau}\cdot)+Y_1(s,\theta_{k\tau}\cdot))(y){\cal D}_rz(s,\theta_{k\tau}\cdot,y)dyds\Big]^2\\
&&+ CE\sum_{i=1}^m\Big[\int^{\infty}_t\int_D e^{\mu_i(t-s)}\phi_i(y)\nabla F^i(s, z(s,\theta_{k\tau}\cdot)+Y_1(s,\theta_{k\tau}\cdot))(y){\cal D}_rz(s,\theta_{k\tau}\cdot,y)dyds\Big]^2\\
&&+ C E\sum_{i=m+1}^\infty\Big[\int_{-\infty}^t\int_D e^{\mu_i(t-s)}\phi_i(y)\nabla F^i(s, z(s,\theta_{k\tau}\cdot)+Y_1(s,\theta_{k\tau}\cdot))(y){\cal D}_rY_1(s,\theta_{k\tau}\cdot,y)dyds\Big]^2\\
&&+ C E\sum_{i=1}^m\Big[\int_{t}^\infty\int_D e^{\mu_i(t-s)}\phi_i(y)\nabla F^i(s, z(s,\theta_{k\tau}\cdot)+Y_1(s,\theta_{k\tau}\cdot))(y){\cal D}_rY_1(s,\theta_{k\tau}\cdot,y)dyds\Big]^2\\
&\leq& C E\sum_{i=m+1}^\infty \Big[\int_{-\infty}^t\int_D e^{\mu_i(t-s)}|\phi_i(y)|^2dyds\cdot\\
&&\hskip 2cm\int_{-\infty}^t\int_D e^{\mu_i(t-s)}|\nabla F^i(s, z(s,\theta_{k\tau}\cdot)+Y_1(s,\theta_{k\tau}\cdot))(y)|^2|{\cal D}_rz(s,\theta_{k\tau}\cdot,y)|^2dyds\Big]\\
 &&+C E\sum_{i=1}^m \Big[\int^{\infty}_t\int_D e^{\mu_i(t-s)}|\phi_i(y)|^2dyds\cdot\\
&&\hskip 2cm\int_t^\infty\int_D e^{2\mu_i(t-s)}|\nabla F^i(s, z(s,\theta_{k\tau}\cdot)+Y_1(s,\theta_{k\tau}\cdot))(y)|^2|{\cal D}_rz(s,\theta_{k\tau}\cdot,y)|^2dyds\Big]\\
&&+ C E\sum_{i=m+1}^\infty \Big[\int_{-\infty}^t\int_D e^{2\mu_i(t-s)}|\phi_i(y)|^2dyds\cdot\\
&&\hskip 2cm\int_{-\infty}^t\int_D |\nabla F^i(s, z(s,\theta_{k\tau}\cdot)+Y_1(s,\theta_{k\tau}\cdot))(y)|^2|{\cal D}_rY_1(s,\theta_{k\tau}\cdot,y)|^2dyds\Big]\\
&&+ C E\sum_{i=1}^m\Big[\int^{\infty}_t\int_D e^{2\mu_i(t-s)}|\phi_i(y)|^2dyds\cdot\\
&&\hskip 1.5cm\int_t^{\infty}\int_D |\nabla F^i(s, z(s,\theta_{k\tau}\cdot)+Y_1(s,\theta_{k\tau}\cdot))(y)|^2|{\cal D}_rY_1(s,\theta_{k\tau}\cdot,y)|^2dyds\Big]\\
&\leq & C(-{1\over {\mu_{m+1}}})||\nabla F||_\infty^2\int_{-\infty}^t\int_D e^{\mu_{m+1}(t-s)}E|{\cal D}_r z(s,\theta_{k\tau}\cdot,y)|^2dyds\\
&&+C{1\over {\mu_{m}}}||\nabla F||_\infty^2\int^{\infty}_t\int_D e^{\mu_{m}(t-s)}E|{\cal D}_r z(s,\theta_{k\tau}\cdot,y)|^2dyds\\
&&+ C(-{2\over {\mu_{m+1}}}) ||\nabla F||_\infty^2\int_{-\infty}^{t+k\tau}\int_DE|{\cal D}_rY_1(s,\cdot,y)|^2dyds\\
&&+ C{2\over {\mu_{m}}} ||\nabla F||_\infty^2\int^{\infty}_{t+k\tau}\int_DE|{\cal D}_rY_1(s,\cdot,y)|^2dyds.
\end{eqnarray*}
Let us first deal with the third and the fourth terms.
When $k=0,1,2,\cdots$, we have $t+k\tau\geq r$ and
\begin{eqnarray*}
&&\int_{-\infty}^{t+k\tau}\int_DE|{\cal D}_rY_1(s,\cdot,y)|^2|dyds\\
&=& \int_{-\infty}^r \int_DE|{\cal D}_rY_1(s,\cdot,y)|^2|dyds+\int_r^{t+k\tau}\int_DE|{\cal D}_rY_1(s,\cdot,y)|^2dyds\\
&=&\int_{-\infty}^r\int_D | \sum_{j=1}^m e^{\mu_j(s-r)}\phi_j(y)\sigma_j(r)|^2dyds+\int_r^{t+k\tau}\int_D | \sum_{j=m+1}^\infty e^{\mu_j(s-r)}\phi_j(y)\sigma_j(r)|^2dyds\\
&=&\int_{-\infty}^r\int_D\sum_{j=1}^m e^{2\mu_j(s-r)}|\phi_j(y)|^2|\sigma_j(r)|^2dyds+\int_r^{t+k\tau}\int_D\sum_{j=m+1}^\infty e^{2\mu_j(s-r)}|\phi_j(y)|^2|\sigma_j(r)|^2dyds\nonumber\\
&\leq&({1\over{2\mu_m}}-{1\over {2\mu_{m+1}}})\sup_{s\in (-\infty,\infty)} \sum_{j=1}^\infty \sigma_j^2(s).
\end{eqnarray*}
When $k=-1,-2,\cdots$, we have $t+k\tau<r$ and 
\begin{eqnarray*}
\int_{-\infty}^{t+k\tau}\int_DE|{\cal D}_rY_1(s,\cdot,y)|^2|dyds
&=&\int_{-\infty}^{t+k\tau}\int_D | \sum_{j=1}^m e^{\mu_j(s-r)}\phi_j(y)\sigma_j(r)|^2dyds\\
&=&\int_{-\infty}^{t+k\tau}\int_D\sum_{j=1}^m e^{2\mu_j(s-r)}|\phi_j(y)|^2|\sigma_j(r)|^2dyds\\
&\leq&{1\over {2\mu_{m}}}\sup_{s\in (-\infty,\infty)} \sum_{j=1}^\infty \sigma_j^2(s).
\end{eqnarray*}
So, \begin{eqnarray*}
\int_{-\infty}^{t+k\tau}\int_DE|{\cal D}_rY_1(s,\cdot,y)|^2|dyds
\leq({1\over{\mu_m}}-{1\over {2\mu_{m+1}}})\sup_{s\in (-\infty,\infty)} \sum_{j=1}^\infty \sigma_j^2(s).
\end{eqnarray*}
Similarly,
\begin{eqnarray*}
\int^{\infty}_{t+k\tau}\int_DE|{\cal D}_rY_1(s,\cdot,y)|^2|dyds
\leq({1\over{2\mu_m}}-{1\over {\mu_{m+1}}})\sup_{s\in (-\infty,\infty)} \sum_{j=1}^\infty \sigma_j^2(s).
\end{eqnarray*}
Therefore, we have
\begin{eqnarray*}
&&E\int_D|{\cal D}_r{\cal M}(z)(t,\theta_{k\tau}\cdot, x)|^2dx\\
&\leq&C(-{1\over {\mu_{m+1}}})||\nabla F||_\infty^2\int_{r-\tau}^r\int_D\sum_{i=0}^\infty e^{\mu_{m+1}(t-s+i\tau)}E|{\cal D}_r z(s-i\tau,\theta_{k\tau}\cdot,y)|^2dyds\\
&&+C(-{1\over {\mu_{m+1}}})||\nabla F||_\infty^2\int_{r}^t\int_D e^{\mu_{m+1}(t-s)}E|{\cal D}_r z(s,\theta_{k\tau}\cdot,y)|^2dyds\\
&&+C{1\over {\mu_{m}}}||\nabla F||_\infty^2\int_t^{r+\tau}\int_D e^{\mu_{m}(t-s)}E|{\cal D}_r z(s,\theta_{k\tau}\cdot,y)|^2dyds\\
&&+C{1\over {\mu_{m}}}||\nabla F||_\infty^2\int_{r+\tau}^{r+2\tau}\int_D\sum_{i=0}^\infty e^{\mu_{m}(t-s-i\tau)}E|{\cal D}_r z(s+i\tau,\theta_{k\tau}\cdot, y)|^2dyds\\
&&+C||\nabla F||_\infty^2\sup_{s\in (-\infty,\infty)} \sum_{j=1}^\infty \sigma_j^2(s)({1\over{\mu_{m+1}^2}}+{1\over {\mu_m^2}})\\
&\leq&C(-{1\over {\mu_{m+1}}})||\nabla F||_\infty^2\sum_{i=0}^\infty e^{\mu_{m+1}i\tau}\int_{r-2\tau}^{r+2\tau}e^{-\beta |t-s|}\int_DE|{\cal D}_r z(s,\theta_{-i\tau+k\tau}\cdot,y)|^2dyds\\
&&+C(-{1\over {\mu_{m+1}}}+{1\over {\mu_{m}}})||\nabla F||_\infty^2\int_{r-2\tau}^{r+2\tau}e^{-\beta |t-s|}E|{\cal D}_r z(s,\theta_{k\tau}\cdot,y)|^2dyds\\
&&+C{1\over {\mu_{m}}}||\nabla F||_\infty^2\sum_{i=0}^\infty e^{-\mu_{m}i\tau}\int_{r-2\tau}^{r+2\tau}e^{-\beta |t-s|}\int_DE|{\cal D}_r z(s,\theta_{i\tau+k\tau}\cdot,y)|^2dyds\\
&&+C||\nabla F||_\infty^2\sup_{s\in (-\infty,\infty)} \sum_{j=1}^\infty \sigma_j^2(s)({1\over{\mu^2_{m+1}}}+{1\over{\mu^2_{m}}})\\
&\leq & A\int_{r-2\tau}^{r+2\tau}e^{-\beta |t-s|}\alpha_r(s)ds+B\\
&=&\alpha_r(t) .
\end{eqnarray*}
Similarly, when $0\leq t< r< \tau$,
\begin{eqnarray*}
&&E\int_D|{\cal D}_r{\cal M}(z)(t,\theta_{k\tau}\cdot,x)|^2dx\\
&\leq&-C{1\over {\mu_{m+1}}}||\nabla F||_\infty^2\int_{r-2\tau}^{r-\tau}\sum_{i=0}^\infty e^{\mu_{m+1}(t-s+i\tau)}\int_DE|{\cal D}_r z(s,s-i\tau,\theta_{k\tau}\cdot,y)|^2dyds\\
&&-C{1\over {\mu_{m+1}}}||\nabla F||_\infty^2\int_{r-\tau}^t e^{\mu_{m+1}(t-s)}\int_DE|{\cal D}_r z(s,\theta_{k\tau}\cdot,y)|^2dyds\\
&&+C{1\over {\mu_{m}}}||\nabla F||_\infty^2\int_t^{r} e^{\mu_{m}(t-s)}E|{\cal D}_r z(s,\theta_{k\tau}\cdot,y)|^2dyds\\
&&+C{1\over {\mu_{m}}}||\nabla F||_\infty^2\int_{r}^{r+\tau}\sum_{i=0}^\infty e^{\mu_{m}(t-s-i\tau)}E|{\cal D}_r z(s,s+i\tau,\theta_{k\tau}\cdot,y)|^2dyds\\
&&+C||\nabla F||_\infty^2\sup_{s\in (-\infty,\infty)} \sum_{j=1}^\infty \sigma_j^2(s)({1\over{\mu^2_{m+1}}}+{1\over{\mu^2_{m}}})\\
&\leq & A\int_{r-2\tau}^{r+2\tau}e^{-\beta |t-s|}\alpha_r(s)ds+B\\
&=&\alpha_r(t) .
\end{eqnarray*}
Therefore, for any $k=0,\pm 1, \pm 2,\cdots$, we have
$$E\int_D|{\cal D}_r{\cal M}(z)(t,\theta_{k\tau}\cdot,x)|^2dx\leq \alpha_r(t).$$
Moreover, the solution $\alpha_r(t)$ of equation (\ref{eqn2.15}) is continuous in $t$, so for $z\in C^0_{\tau,\alpha}((-\infty,+\infty),L^2(D,{\cal D}^{1,2}))$, there exists a constant $\alpha_1$ such that for any $t, r \in [0, \tau), k=0,\pm 1, \pm 2,\cdots$,
\begin{eqnarray*}
E\int_D|{\cal D}_r z(t,\theta_{k\tau}\cdot, x)|^2dx \leq \alpha_1, \ and\ \ 
E\int_D|{\cal D}_r {\cal M}(z)(t,\theta_{k\tau}\cdot, x)|^2 dx\leq \alpha_1, 
\end{eqnarray*}
 Now suppose there exists $L_2\geq 0$ such that  for any $r_1, r_2, s\in [0,\tau),  k=0,\pm 1, \pm 2,\cdots$,
\begin{eqnarray*}
{1\over {|r_1-r_2|}}\int_D{{E|{\cal D}_{r_1} z(s,\theta_{k\tau}\cdot, x)-{\cal D}_{r_2} z(s,\theta_{k\tau}\cdot, x)|^2}}dx\leq L_2.
\end{eqnarray*}
Then we have when $0\leq r_1<r_2\leq t<\tau$,  $k=0,\pm 1, \pm 2,\cdots$
\begin{eqnarray*}
&&{1\over {|r_1-r_2|}}\int_DE|{\cal D}_{r_1}{\cal M}(z)(t,\theta_{k\tau}\cdot, x)-{\cal D}_{r_2}{\cal M}(z)(t,\theta_{k\tau}\cdot, x)|^2dx\\
&\leq &{C\over {|r_1-r_2|}}\int_D\Big\{E\big|\int_{-\infty}^t \int_D\sum_{i=m+1}^\infty e^{\mu_i(t-s)}\phi_i(x)\phi_i(y)\nabla F^i(s, z(s,\theta_{k\tau}\cdot)+Y_1(s,\theta_{k\tau}\cdot))(y)\\
&&\hskip 4cm \cdot({\cal D}_{r_1} z(s,\theta_{k\tau}\cdot, y)-{\cal D}_{r_2} z(s,\theta_{k\tau}\cdot, y))dyds\big|^2\\
&&+E\big|\int^{\infty}_t \int_D\sum_{i=1}^m e^{\mu_i(t-s)}\phi_i(x)\phi_i(y)\nabla F^i(s, z(s,\theta_{k\tau}\cdot)+Y_1(s,\theta_{k\tau}\cdot))(y)\\
&&\hskip 4cm\cdot({\cal D}_{r_1} z(s,\theta_{k\tau}\cdot, y)-{\cal D}_{r_2} z(s,\theta_{k\tau}\cdot, y))dyds\big|^2\\
&&+E\big|\int_{r_1}^t \int_D\sum_{i=m+1}^\infty e^{\mu_i(t-s)}\phi_i(x)\phi_i(y)\nabla F^i(s, z(s,\theta_{k\tau}\cdot)+Y_1(s,\theta_{k\tau}\cdot))(y){\cal D}_{r_1}Y_1(s,\theta_{k\tau}\cdot, y)dyds\\
&&\hskip 0.5cm-\int_{r_2}^t \int_D\sum_{i=m+1}^\infty e^{\mu_i(t-s)}\phi_i(x)\phi_i(y)\nabla F^i(s, z(s,\theta_{k\tau}\cdot)+Y_1(s,\theta_{k\tau}\cdot))(y){\cal D}_{r_2}Y_1(s,\theta_{k\tau}\cdot, y)dyds\big|^2\\
&&+E\big|\int^{r_1}_{-\infty} \int_D\sum_{i=m+1}^\infty e^{\mu_i(t-s)}\phi_i(x)\phi_i(y)\nabla F^i(s, z(s,\theta_{k\tau}\cdot)+Y_1(s,\theta_{k\tau}\cdot))(y){\cal D}_{r_1}Y_1(s,\theta_{k\tau}\cdot, y)dyds\\
&&\hskip 0.5cm-\int_{-\infty}^{r_2} \int_D\sum_{i=m+1}^\infty e^{\mu_i(t-s)}\phi_i(x)\phi_i(y)\nabla F^i(s, z(s,\theta_{k\tau}\cdot)+Y_1(s,\theta_{k\tau}\cdot))(y){\cal D}_{r_2}Y_1(s,\theta_{k\tau}\cdot, y))dyds\big|^2\\
&&+E\big|\int_{t}^{\infty} \sum_{i=1}^m e^{\mu_i(t-s)}\phi_i(x)\phi_i(y)\nabla F^i(s, z(s,\theta_{k\tau}\cdot)+Y_1(s,\theta_{k\tau}\cdot))(y)\\
&&\hskip 4cm\cdot ({\cal D}_{r_1}Y_1(s,\theta_{k\tau}\cdot, y)-{\cal D}_{r_2}Y_1(s,\theta_{k\tau}\cdot, y))dyds\big|^2\Big\}dx\\
&:=&A_1+A_2+A_3+A_4+A_5.
\end{eqnarray*}
We will estimate them in the following. We first have that
\begin{eqnarray*}
A_1&\leq& {C\over {|r_1-r_2|}}\int_D E\Big|\int_{-\infty}^t\int_D\sum_{i=m+1}^\infty e^{\mu_i(t-s)}\phi_i(x)\phi_i(y)\nabla F^i(s, z(s,\theta_{k\tau}\cdot)+Y_1(s,\theta_{k\tau}\cdot))(y)\\
&&\hskip 4cm\cdot ({\cal D}_{r_1} z(s,\theta_{k\tau}\cdot, y)-{\cal D}_{r_2} z(s,\theta_{k\tau}\cdot, y))dyds\Big|^2dx\\
&\leq&  {C\over {|r_1-r_2|}}E\sum_{i=m+1}^\infty\Big[\int_{-\infty}^t\int_D e^{\mu_i(t-s)}|\phi_i(y)|^2dyds\\
&&\cdot\int_{-\infty}^t\int_D e^{\mu_i(t-s)}|\nabla F^i(s, z(s,\theta_{k\tau}\cdot)+Y_1(s,\theta_{k\tau}\cdot))(y)|^2 
 |{\cal D}_{r_1} z(s,\theta_{k\tau}\cdot, y)-{\cal D}_{r_2} z(s,\theta_{k\tau}\cdot, y)|^2dyds\Big]\\
 &\leq& {C\over {|r_1-r_2|}}||\nabla F||_{\infty}^2 (-{1\over {\mu_{m+1}}})\int_{-\infty}^t\int_D e^{\mu_{m+1}(t-s)}E|{\cal D}_{r_1} z(s,\theta_{k\tau}\cdot, y)-{\cal D}_{r_2} z(s,\theta_{k\tau}\cdot, y)|^2dyds \\
 &\leq& {C\over {|r_1-r_2|}}||\nabla F||_{\infty}^2 (-{1\over {\mu_{m+1}}})\Big[\int_0^t \int_D e^{\mu_{m+1}(t-s)}E|{\cal D}_{r_1} z(s,\theta_{+k\tau}\cdot, y)-{\cal D}_{r_2} z(s,\theta_{k\tau}\cdot, y)|^2dyds\\ 
 &&+\int_0^\tau\int_D \sum_{i=0}^\infty e^{\mu_{m+1}(t-s+\tau+i\tau)}E|{\cal D}_{r_1} z(s,\theta_{-(i+1)\tau+k\tau}\cdot, y)-{\cal D}_{r_2} z(s,\theta_{-(i+1)\tau+k\tau}\cdot, y)|^2dyds\Big]\\
&\leq&{C\over {\mu^2_{m+1}}}||\nabla F||_{\infty}^2L_2[1+\sum_{i=0}^\infty e^{\mu_{m+1}i\tau}].
\end{eqnarray*}
Similarly,
\begin{eqnarray*}
A_2\leq{C\over {\mu^2_{m}}}||\nabla F||_{\infty}^2L_2[1+\sum_{i=0}^\infty e^{-\mu_{m}i\tau}].
\end{eqnarray*}
For $A_3$, using  Cauchy-Schwarz inequality again, we have
\begin{eqnarray*}
A_3&=&{C\over {|r_1-r_2|}}\int_DE\Big|\int_{r_1}^{r_2}\int_D\sum_{i=m+1}^\infty e^{\mu_i(t-s)}\phi_i(x)\phi_i(y)\\
&&\hskip 3cm \cdot\nabla F^i(s, z(s,\theta_{k\tau}\cdot)+Y_1(s,\theta_{k\tau}\cdot))(y)D_{r_1}Y_1(s,\theta_{k\tau}\cdot,y)dyds\Big|^2dx\\
&&+{C\over {|r_1-r_2|}}\int_DE\Big|\int_{r_2}^t\int_D\sum_{i=m+1}^\infty e^{\mu_i(t-s)}\phi_i(x)\phi_i(y)\nabla F^i(s, z(s)+Y_1(s))(y)\\
&&\hskip 3cm \cdot (D_{r_1}Y_1(s,\theta_{k\tau}\cdot,y)-D_{r_2}Y_1(s,\theta_{k\tau}\cdot,y)dyds\Big|^2dx\\
&\leq& {C\over {|r_1-r_2|}}\sum_{i=m+1}^\infty E\int_{r_1}^{r_2}\int_D e^{2\mu_i(t-s)}|\phi_i(y)|^2dyds\\
&&\hskip 2cm\cdot\int_{r_1}^{r_2}\int_D |\nabla F^i(s, z(s)+Y_1(s))(y)|^2\Big|D_{r_1}Y_1(s,\theta_{k\tau}\cdot,y)\Big|^2dyds\\
&&+{C\over {|r_1-r_2|}}\sum_{i=m+1}^\infty E\int^t_{r_2}\int_D e^{\mu_i(t-s)}|\phi_i(y)|^2dyds\\
&&\hskip1cm \cdot\int^t_{r_2}\int_D e^{\mu_i(t-s)} |\nabla F^i(s, z(s)+Y_1(s))(y)|^2 |D_{r_1}Y_1(s,\theta_{k\tau}\cdot,y)-D_{r_2}Y_1(s,\theta_{k\tau}\cdot,y)|^2dyds\\
&\leq&{C\over {|r_1-r_2|}}(r_2-r_1)||\nabla F||_\infty^2\int_{r_1}^{r_2}\int_D E|D_{r_1}Y_1(s,\theta_{k\tau}\cdot,y)|^2dyds\\
&&+ {C\over {|r_1-r_2|}}(t-r_2)||\nabla F||_\infty^2\cdot\int_{r_2}^t\int_D E| D_{r_1}Y_1(s,\theta_{k\tau}\cdot,y)-D_{r_2}Y_1(s,\theta_{k\tau}\cdot,y)|^2dyds\\
%
\end{eqnarray*}
Note that
\begin{eqnarray*}
\int_{r_1}^{r_2}\int_D E|D_{r_1}Y_1(s,\theta_{k\tau}\cdot,y)|^2dyds=\int_{r_1+k\tau}^{r_2+k\tau}\int_D E|D_{r_1}Y_1(s,\cdot,y)|^2dyds,
\end{eqnarray*}
so when $k=0,1, 2,\cdots$, we have
\begin{eqnarray*}
\int_{r_1+k\tau}^{r_2+k\tau}\int_D E|D_{r_1}Y_1(s,\cdot,y)|^2dyds
&=&\int_{r_1+k\tau}^{r_2+k\tau}\int_D |\sum_{j=m+1}^\infty e^{\mu_j(s-r_1)} \phi_j(y)\sigma_j(r_1)|^2dyds\\
&=&\int_{r_1+k\tau}^{r_2+k\tau}\int_D \sum_{j=m+1}^\infty e^{2\mu_j(s-r_1)} |\phi_j(y)|^2|\sigma_j(r_1)|^2dyds\\
&\leq& |r_2-r_1| \sup_{s\in (-\infty,\infty)} \sum_{j=1}^\infty |\sigma_j(s)|^2.
\end{eqnarray*}
When $k=-1,-2,\cdots$,  we have  $r_2+k\tau< r_1$ and 
\begin{eqnarray*}
\int_{r_1+k\tau}^{r_2+k\tau}\int_D E|D_{r_1}Y_1(s,\cdot,y)|^2dyds
&=&\int_{r_1+k\tau}^{r_2+k\tau}\int_D \sum_{j=1}^m e^{2\mu_j(s-r_1)} |\phi_j(y)|^2|\sigma_j(r_1)|^2dyds\\
&\leq& |r_2-r_1|\sup_{s\in (-\infty,\infty)} \sum_{j=1}^\infty |\sigma_j(s)|^2.
\end{eqnarray*}
Therefore 
\begin{eqnarray}\label{eqn2.12}
\int_{r_1}^{r_2}\int_D E|D_{r_1}Y_1(s,\theta_{k\tau}\cdot,y)|^2dyds \leq  |r_2-r_1|\sup_{s\in (-\infty,\infty)} \sum_{j=1}^\infty |\sigma_j(s)|^2.
\end{eqnarray}
Similarly,
\begin{eqnarray*}
&&\int_{r_2}^t\int_D E| D_{r_1}Y_1(s,\theta_{k\tau}\cdot,y)-D_{r_2}Y_1(s,\theta_{k\tau}\cdot,y)|^2dyds\\
&=&\int_{r_2+k\tau}^{t+k\tau}\int_D E| D_{r_1}Y_1(s,\cdot,y)-D_{r_2}Y_1(s,\cdot,y)|^2dyds
\end{eqnarray*}
When $k=0, 1, 2, \cdots$, we have
\begin{eqnarray*}
&&\int_{r_2+k\tau}^{t+k\tau}\int_D E| D_{r_1}Y_1(s,\cdot,y)-D_{r_2}Y_1(s,\cdot,y)|^2dyds\\
&=&\int_{r_2+k\tau}^{t+k\tau}\int_D \Big|\sum_{j=m+1}^\infty \big(e^{\mu_j(s-r_1)}\sigma_j(r_1)-e^{\mu_j(s-r_1)}\sigma_j(r_2)+e^{\mu_j(s-r_1)}\sigma_j(r_2)-e^{\mu_j(s-r_2)}\sigma_j(r_2)\big)\phi_j(y)\Big|^2dyds\\
&=&\int_{r_2+k\tau}^{t+k\tau}\int_D \sum_{j=m+1}^\infty \big|e^{\mu_j(s-r_1)}(\sigma_j(r_1)-\sigma_j(r_2))+(e^{\mu_j(s-r_1)}-e^{\mu_j(s-r_2)})\sigma_j(r_2)\big|^2|\phi_j(y)|^2dyds\\
&\leq&2\int_{r_2+k\tau}^{t+k\tau}\int_D\sum_{j=m+1}^\infty |\sigma_j(r_1)-\sigma_j(r_2)|^2|\phi_j(y)|^2dyds\\
&&+2\int_{r_2+k\tau}^{t+k\tau}\int_D\sum_{j=m+1}^\infty |e^{\mu_j(s-r_1)}-e^{\mu_j(s-r_2)}|^2|\phi_j(y)|^2|\sigma_j(r_2)|^2dyds\\
&\leq&2L_1|r_2-r_1|(t-r_2)+|r_2-r_1|(t-r_2) \sup_{s\in (-\infty,\infty)} \sum_{j=1}^\infty |\sigma_j(s)|^2.
\end{eqnarray*}
When $k=-1, -2, \cdots$, we have $r_2+k\tau< t+k\tau<r_1<r_2\leq t<\tau$ and 
\begin{eqnarray*}
&&\int_{r_2+k\tau}^{t+k\tau}\int_D E| D_{r_1}Y_1(s,\cdot,y)-D_{r_2}Y_1(s,\cdot,y)|^2dyds\\
&=&\int_{r_2+k\tau}^{t+k\tau}\int_D \Big|\sum_{j=1}^m \big(e^{\mu_j(s-r_1)}\sigma_j(r_1)-e^{\mu_j(s-r_1)}\sigma_j(r_2)+e^{\mu_j(s-r_1)}\sigma_j(r_2)-e^{\mu_j(s-r_2)}\sigma_j(r_2)\big)\phi_j(y)\Big|^2dyds\\
&=&\int_{r_2+k\tau}^{t+k\tau}\int_D \sum_{j=1}^m \big|e^{\mu_j(s-r_1)}(\sigma_j(r_1)-\sigma_j(r_2))+(e^{\mu_j(s-r_1)}-e^{\mu_j(s-r_2)})\sigma_j(r_2)\big|^2|\phi_j(y)|^2dyds\\
&\leq&2\int_{r_2+k\tau}^{t+k\tau}\int_D\sum_{j=1}^m |\sigma_j(r_1)-\sigma_j(r_2)|^2|\phi_j(y)|^2dyds\\
&&+2\int_{r_2+k\tau}^{t+k\tau}\int_D\sum_{j=1}^m |e^{\mu_j(s-r_1)}-e^{\mu_j(s-r_2)}|^2|\phi_j(y)|^2|\sigma_j(r_2)|^2dyds\\
&\leq&2L_1|r_2-r_1|(t-r_2)+|r_2-r_1| \sup_{s\in (-\infty,\infty)} \sum_{j=1}^\infty |\sigma_j(s)|^2.
\end{eqnarray*}
Therefore,
\begin{eqnarray}\label{eqn2.13}, 
&&\int_{r_2}^t\int_D E| D_{r_1}Y_1(s,\theta_{k\tau}\cdot,y)-D_{r_2}Y_1(s,\theta_{k\tau}\cdot,y)|^2dyds\nonumber\\
&\leq&2L_1|r_2-r_1|(t-r_2)+|r_2-r_1| \sup_{s\in (-\infty,\infty)} \sum_{j=1}^\infty |\sigma_j(s)|^2
\end{eqnarray}
With the estimate (\ref{eqn2.12}) and (\ref{eqn2.13}), we have
\begin{eqnarray*}
A_3
&\leq&{C\over {|r_1-r_2|}}||\nabla F||_{\infty}^2(r_2-r_1)^2\sup_{s\in(-\infty,\infty)}  \sum_{j=1}^\infty|\sigma_j(s)|^2\\
&&+{{C}\over {|r_1-r_2|}}\Big[||\nabla F||_{\infty}^2(t-r_2)^22L_1|r_2-r_1|+||\nabla F||_{\infty}^2|t-r_2||r_2-r_1|\sup_{s\in (-\infty,\infty)} \sum_{j=1}^\infty \sigma_j^2(s)\Big]\\
&\leq & C||\nabla F||_{\infty}^2\tau\sup_{s\in(-\infty,\infty)}  \sum_{j=1}^\infty|\sigma_j(s)|^2+C||\nabla F||_{\infty}^2(2L_1\tau^2+\tau\sup_{s\in(-\infty,\infty)}  \sum_{j=1}^\infty|\sigma_j(s)|^2)\\
&<&\infty.
\end{eqnarray*}
About $A_4$,
\begin{eqnarray*}
A_4
&=&{C\over {|r_1-r_2|}}\int_DE\Big|\int_{-\infty}^{r_1}\int_D\sum_{i=m+1}^\infty e^{\mu_i(t-s)}\phi_i(x)\phi_i(y)\nabla F^i(s, z(s)+Y_1(s))(y)\\
&&\hskip 3cm \cdot (D_{r_1}Y_1(s,\theta_{k\tau}\cdot,y)-D_{r_2}Y_1(s,\theta_{k\tau}\cdot,y)dyds\Big|^2dx\\
&&+{C\over {|r_1-r_2|}}\int_DE\Big|\int_{r_1}^{r_2}\int_D\sum_{i=m+1}^\infty e^{\mu_i(t-s)}\phi_i(x)\phi_i(y)\\
&&\hskip 3cm \cdot\nabla F^i(s, z(s,\theta_{k\tau}\cdot)+Y_1(s,\theta_{k\tau}\cdot))(y)D_{r_2}Y_1(s,\theta_{k\tau}\cdot,y)dyds\Big|^2dx\\
&\leq&{C\over {|r_1-r_2|}}\sum_{i=m+1}^\infty E\int_{-\infty}^{r_1}\int_D e^{2\mu_i(t-s)}|\phi_i(y)|^2dyds\\
&&\hskip1cm \cdot\int_{-\infty}^{r_1}\int_D  |\nabla F^i(s, z(s)+Y_1(s))(y)|^2 |D_{r_1}Y_1(s,\theta_{k\tau}\cdot,y)-D_{r_2}Y_1(s,\theta_{k\tau}\cdot,y)|^2dyds\\
&&+ {C\over {|r_1-r_2|}}\sum_{i=m+1}^\infty E\int_{r_1}^{r_2}\int_D e^{2\mu_i(t-s)}|\phi_i(y)|^2dyds\\
&&\hskip 2cm\cdot\int_{r_1}^{r_2}\int_D |\nabla F^i(s, z(s)+Y_1(s))(y)|^2\Big|D_{r_2}Y_1(s,\theta_{k\tau}\cdot,y)\Big|^2dyds\\
&\leq& {C\over {|r_1-r_2|}}(-{1\over{2\mu_{m+1}}})||\nabla F||_\infty^2\int_{-\infty}^{r_1}\int_D E| D_{r_1}Y_1(s,\theta_{k\tau}\cdot,y)-D_{r_2}Y_1(s,\theta_{k\tau}\cdot,y)|^2dyds\\
&&+{C\over {|r_1-r_2|}}(r_2-r_1)||\nabla F||_\infty^2\int_{r_1}^{r_2}\int_D E|D_{r_2}Y_1(s,\theta_{k\tau}\cdot,y)|^2dyds.
\end{eqnarray*}
Similar to (\ref{eqn2.12}), 
\begin{eqnarray}\label{2.14}
\int_{r_1}^{r_2}\int_D E|D_{r_2}Y_1(s,\theta_{k\tau}\cdot,y)|^2dyds\leq  |r_2-r_1|\sup_{s\in (-\infty,\infty)} \sum_{j=1}^\infty |\sigma_j(s)|^2.
\end{eqnarray}
Secondly, 
\begin{eqnarray*}
 &&\int_{-\infty}^{r_1}\int_D E| D_{r_1}Y_1(s,\theta_{k\tau}\cdot,y)-D_{r_2}Y_1(s,\theta_{k\tau}\cdot,y)|^2dyds\\
 &=&\int_{-\infty}^{r_1+k\tau}\int_D E| D_{r_1}Y_1(s,\cdot,y)-D_{r_2}Y_1(s,\cdot,y)|^2dyds.
\end{eqnarray*}
When $k=0, -1,-2,\cdots$,  we have
\begin{eqnarray*}
 &&\int_{-\infty}^{r_1+k\tau}\int_D E| D_{r_1}Y_1(s,\cdot,y)-D_{r_2}Y_1(s,\cdot,y)|^2dyds\\
 &=&\int_{-\infty}^{r_1+k\tau}\int_D\Big|\sum_{j=1}^m \big(e^{\mu_j(s-r_1)}\sigma_j(r_1)-e^{\mu_j(s-r_1)}\sigma_j(r_2)+e^{\mu_j(s-r_1)}\sigma_j(r_2)-e^{\mu_j(s-r_2)}\sigma_j(r_2)\big)\phi_j(y)\Big|^2dyds\\
&=&\int_{-\infty}^{r_1+k\tau}\int_D \sum_{j=1}^m \big|e^{\mu_j(s-r_1)}(\sigma_j(r_1)-\sigma_j(r_2))+(e^{\mu_j(s-r_1)}-e^{\mu_j(s-r_2)})\sigma_j(r_2)\big|^2|\phi_j(y)|^2dyds\\
&\leq&2\int_{-\infty}^{r_1+k\tau}\int_D\sum_{j=1}^m e^{\mu_{m}(s-r_1)}|\sigma_j(r_1)-\sigma_j(r_2)|^2|\phi_j(y)|^2dyds\\
&&+2\int_{-\infty}^{r_1+k\tau}\int_D\sum_{j=1}^m |e^{\mu_j(s-r_1)}-e^{\mu_j(s-r_2)}|^2|\phi_j(y)|^2|\sigma_j(r_2)|^2dyds\\
&\leq&2L_1{1\over{\mu_m}}|r_2-r_1|+|r_2-r_1| \sup_{s\in (-\infty,\infty)} \sum_{j=1}^\infty |\sigma_j(s)|^2.
\end{eqnarray*}
When $k=1, 2, \cdots$, we have  $r_1+k\tau> r_2$ and 
\begin{eqnarray*}
 &&\int_{-\infty}^{r_1+k\tau}\int_D E| D_{r_1}Y_1(s,\cdot,y)-D_{r_2}Y_1(s,\cdot,y)|^2dyds\\
 &=&\int_{-\infty}^{r_1}\int_D E| D_{r_1}Y_1(s,\cdot,y)-D_{r_2}Y_1(s,\cdot,y)|^2dyds+\int_{r_1}^{r_2}\int_D E| D_{r_1}Y_1(s,\cdot,y)-D_{r_2}Y_1(s,\cdot,y)|^2dyds\\
 &&+\int_{r_2}^{r_1+k\tau}\int_D E| D_{r_1}Y_1(s,\cdot,y)-D_{r_2}Y_1(s,\cdot,y)|^2dyds\
\end{eqnarray*}
Let us estimate them separately. About the first term,
\begin{eqnarray*}
&&\int_{-\infty}^{r_1}\int_D E| D_{r_1}Y_1(s,\cdot,y)-D_{r_2}Y_1(s,\cdot,y)|^2dyds\\
&=&\int_{-\infty}^{r_1}\int_D\Big|\sum_{j=1}^m \big(e^{\mu_j(s-r_1)}\sigma_j(r_1)-e^{\mu_j(s-r_1)}\sigma_j(r_2)+e^{\mu_j(s-r_1)}\sigma_j(r_2)-e^{\mu_j(s-r_2)}\sigma_j(r_2)\big)\phi_j(y)\Big|^2dyds\\
&=&\int_{-\infty}^{r_1}\int_D \sum_{j=1}^m \big|e^{\mu_j(s-r_1)}(\sigma_j(r_1)-\sigma_j(r_2))+(e^{\mu_j(s-r_1)}-e^{\mu_j(s-r_2)})\sigma_j(r_2)\big|^2|\phi_j(y)|^2dyds\\
&\leq&2\int_{-\infty}^{r_1}\int_D\sum_{j=1}^m e^{\mu_{m}(s-r_1)}|\sigma_j(r_1)-\sigma_j(r_2)|^2|\phi_j(y)|^2dyds\\
&&+2\int_{-\infty}^{r_1}\int_D\sum_{j=1}^m |e^{\mu_j(s-r_1)}-e^{\mu_j(s-r_2)}|^2|\phi_j(y)|^2|\sigma_j(r_2)|^2dyds\\
&\leq&2L_1{1\over{\mu_m}}|r_2-r_1|+|r_2-r_1| \sup_{s\in (-\infty,\infty)} \sum_{j=1}^\infty |\sigma_j(s)|^2.
\end{eqnarray*}
About the second term,
\begin{eqnarray*}
&&\int_{r_1}^{r_2}\int_D E| D_{r_1}Y_1(s,\cdot,y)-D_{r_2}Y_1(s,\cdot,y)|^2dyds\\
&\leq &2\int_{r_1}^{r_2}\int_D E| D_{r_1}Y_1(s,\cdot,y)|^2dyds+2\int_{r_1}^{r_2}\int_D E|D_{r_2}Y_1(s,\cdot,y)|^2dyds\\
&\leq& 2 |r_2-r_1|\sup_{s\in (-\infty,\infty)} \sum_{j=1}^\infty |\sigma_j(s)|^2
\end{eqnarray*}
About the third term,
\begin{eqnarray*}
&&\int_{r_2}^{r_1+k\tau}\int_D E| D_{r_1}Y_1(s,\cdot,y)-D_{r_2}Y_1(s,\cdot,y)|^2dyds\\
&=&\int_{-\infty}^{r_1+k\tau}\int_D\Big|\sum_{j=m+1}^\infty \big(e^{\mu_j(s-r_1)}\sigma_j(r_1)-e^{\mu_j(s-r_1)}\sigma_j(r_2)+e^{\mu_j(s-r_1)}\sigma_j(r_2)-e^{\mu_j(s-r_2)}\sigma_j(r_2)\big)\phi_j(y)\Big|^2dyds\\
&=&\int_{-\infty}^{r_1+k\tau}\int_D \sum_{j=m+1}^\infty \big|e^{\mu_j(s-r_1)}(\sigma_j(r_1)-\sigma_j(r_2))+(e^{\mu_j(s-r_1)}-e^{\mu_j(s-r_2)})\sigma_j(r_2)\big|^2|\phi_j(y)|^2dyds\\
&\leq&2\int_{-\infty}^{r_1+k\tau}\int_D\sum_{j=m+1}^\infty e^{\mu_{m}(s-r_1)}|\sigma_j(r_1)-\sigma_j(r_2)|^2|\phi_j(y)|^2dyds\\
&&+2\int_{-\infty}^{r_1+k\tau}\int_D\sum_{j=m+1}^\infty |e^{\mu_j(s-r_1)}-e^{\mu_j(s-r_2)}|^2|\phi_j(y)|^2|\sigma_j(r_2)|^2dyds\\
&\leq&2L_1(-{1\over{\mu_{m+1}}})|r_2-r_1|+|r_2-r_1| \sup_{s\in (-\infty,\infty)} \sum_{j=1}^\infty |\sigma_j(s)|^2.
\end{eqnarray*}
Therefore,
\begin{eqnarray}\label{2.15}
 &&\int_{-\infty}^{r_1}\int_D E| D_{r_1}Y_1(s,\theta_{k\tau}\cdot,y)-D_{r_2}Y_1(s,\theta_{k\tau}\cdot,y)|^2dyds\nonumber\\
 &\leq&2L_1({1\over {\mu_m}}-{1\over{\mu_{m+1}}})|r_2-r_1|+|r_2-r_1| \sup_{s\in (-\infty,\infty)} \sum_{j=1}^\infty |\sigma_j(s)|^2.
 \end{eqnarray}
 With (\ref{2.14}) and ({\ref{2.15}), we have
\begin{eqnarray*}
A_4&\leq&{C\over {|r_1-r_2|}}||\nabla F||_\infty^2\Big[(-{1\over {2\mu_{m+1}}}) (2L_1({1\over {\mu_m}}-{1\over{\mu_{m+1}}})|r_2-r_1| +|r_2-r_1| \sup_{s\in(-\infty,\infty)}  \sum_{j=1}^\infty|\sigma_j(s)|^2)\\
&&\hskip 2cm+|r_2-r_1|^2 \sup_{s\in(-\infty,\infty)}  \sum_{j=1}^\infty|\sigma_j(s)|^2\Big]\\
&\leq& C||\nabla F||_\infty^2\Big(({1\over {\mu_{m+1}^2}}-{1\over {\mu_{m+1}\mu_m}})L_1+(-{1\over {2\mu_{m+1}}})  \sup_{s\in(-\infty,\infty)}  \sum_{j=1}^\infty|\sigma_j(s)|^2+ \tau\cdot\sup_{s\in(-\infty,\infty)}  \sum_{j=1}^\infty|\sigma_j(s)|^2\Big)\\
&<&\infty.
\end{eqnarray*}
As for $A_5$, similarly to $A_4$, we have
\begin{eqnarray*}
A_5
\leq C||\nabla F||_\infty^2\Big(({1\over {\mu_{m}^2}}-{1\over {\mu_{m+1}\mu_m}})L_1+{1\over {2\mu_{m}}} \sup_{s\in(-\infty,\infty)}  \sum_{j=1}^\infty|\sigma_j(s)|^2\Big)
<\infty.
\end{eqnarray*}
So, when $0\leq r_1<r_2\leq t<\tau$,  
\begin{eqnarray*}
{1\over {|r_1-r_2|}}\int_DE|{\cal D}_{r_1}{\cal M}(z)(t,\theta_{k\tau}\cdot, x)-{\cal D}_{r_2}{\cal M}(z)(t,\theta_{k\tau}\cdot, x)|^2dx\leq \hat C.
\end{eqnarray*}
When $0\leq r_1<t<r_2<\tau$, 
\begin{eqnarray*}
&&{1\over {|r_1-r_2|}}\int_DE|{\cal D}_{r_2}{\cal M}(z)(t,\theta_{k\tau}\cdot,x)-{\cal D}_{r_1}{\cal M}(z)(t,\theta_{k\tau}\cdot,x)|^2dx\\
&\leq &{C\over {|r_1-r_2|}}\int_D\Big\{E\big|\int_{-\infty}^t \int_D\sum_{i=m+1}^\infty e^{\mu_i(t-s)}\phi_i(x)\phi_i(y)\cdot\nabla F^i(s, z(s,\theta_{k\tau}\cdot)+Y_1(s,\theta_{k\tau}\cdot))(y)\\
&&\hskip2cm({\cal D}_{r_2} z(s,\theta_{k\tau}\cdot,y)-{\cal D}_{r_1} z(s,\theta_{k\tau}\cdot,y))dyds\big|^2\\
&&+E\big|\int^{\infty}_t \int_D\sum_{i=1}^m e^{\mu_i(t-s)}\phi_i(x)\phi_i(y)\nabla F^i(s, z(s,\theta_{k\tau}\cdot)+Y_1(s,\theta_{k\tau}\cdot))(y)\\
&&\hskip 2cm \cdot({\cal D}_{r_2} z(s,\theta_{k\tau}\cdot,y)-{\cal D}_{r_1} z(s,\theta_{k\tau}\cdot,y))dyds\big|^2\\
&&+E\big|-\int_{-\infty}^t \int_D\sum_{i=m+1}^\infty e^{\mu_i(t-s)}\phi_i(x)\phi_i(y)\nabla F^i(s, z(s,\theta_{k\tau}\cdot)+Y_1(s,\theta_{k\tau}\cdot))(y){\cal D}_{r_2}Y_1(s,\theta_{k\tau}\cdot,y)dyds\\
&&\hskip 1cm+\int_{-\infty}^{r_1} \int_D\sum_{i=m+1}^\infty e^{\mu_i(t-s)}\phi_i(x)\phi_i(y)\nabla F^i(s, z(s,\theta_{k\tau}\cdot)+Y_1(s,\theta_{k\tau}\cdot))(y){\cal D}_{r_1}Y_1(s,\theta_{k\tau}\cdot,y)dyds\big|^2\\
&&+E\big|\int^{r_1}_{-\infty} \int_D\sum_{i=m+1}^\infty e^{\mu_i(t-s)}\phi_i(x)\phi_i(y)\nabla F^i(s, z(s,\theta_{k\tau}\cdot)+Y_1(s,\theta_{k\tau}\cdot))(y){\cal D}_{r_1}Y_1(s,\theta_{k\tau}\cdot,y)dyds\\
&&\hskip 1cm-\int_{-\infty}^{r_2} \int_D\sum_{i=m+1}^\infty e^{\mu_i(t-s)}\phi_i(x)\phi_i(y)\nabla F^i(s, z(s,\theta_{k\tau}\cdot)+Y_1(s,\theta_{k\tau}\cdot))(y){\cal D}_{r_2}Y_1(s,\theta_{k\tau}\cdot,y)dyds\big|^2\\
&&+E\big|\int_t^{r_2} \int_D\sum_{i=1}^m e^{\mu_i(t-s)}\phi_i(x)\phi_i(y)\nabla F^i(s, z(s,y)+Y_1(s,\theta_{k\tau}\cdot,y)){\cal D}_{r_2}Y_1(s,\theta_{k\tau}\cdot,y)dyds\big|^2\\
&&+E\big|\int_{r_1}^t \int_D\sum_{i=m+1}^\infty e^{\mu_i(t-s)}\phi_i(x)\phi_i(y)\nabla F^i(s, z(s,\theta_{k\tau}\cdot)+Y_1(s,\theta_{k\tau}\cdot))(y){\cal D}_{r_1}Y_1(s,\theta_{k\tau}\cdot,y)dyds\big|^2\\
&&+E\big|-\int_{r_2}^{\infty} \int_D\sum_{i=1}^m e^{\mu_i(t-s)}\phi_i(x)\phi_i(y)\nabla F^i(s, z(s,\theta_{k\tau}\cdot)+Y_1(s,\theta_{k\tau}\cdot))(y){\cal D}_{r_2}Y_1(s,\theta_{k\tau}\cdot,y)dyds\\
&&\hskip 1cm+\int_t^{\infty} \int_D\sum_{i=1}^m e^{\mu_i(t-s)}\phi_i(x)\phi_i(y)\nabla F^i(s, z(s,\theta_{k\tau}\cdot)+Y_1(s,\theta_{k\tau}\cdot))(y){\cal D}_{r_1}Y_1(s,\theta_{k\tau}\cdot,y)dyds\big|^2\Big\}dx
\end{eqnarray*}
Thus using a similar method as before, we can see that
\begin{eqnarray*}
&&{1\over {|r_1-r_2|}}\int_DE|{\cal D}_{r_2}{\cal M}(z)(t,\theta_{k\tau}\cdot,x)-{\cal D}_{r_1}{\cal M}(z)(t,\theta_{k\tau}\cdot,x)|^2dx\\
&\leq& C||\nabla F||_\infty^2\Big\{({1\over {\mu^2_{m+1}}}+{1\over {\mu^2_{m}}})L_2(1+\sum_{i=0}^\infty e^{\mu_{m+1}i\tau}+\sum_{i=0}^\infty e^{-\mu_{m}i\tau})+({1\over{\mu^2_{m+1}}}+{1\over{\mu^2_{m}}}-{1\over{\mu_m\mu_{m+1}}}+2\tau^2)L_1\\
&&+(-{1\over{2\mu_{m+1}}}+{1\over{2\mu_{m}}}+4\tau)\sup_{s\in(-\infty,\infty)}  \sum_{j=1}^\infty|\sigma_j(s)|^2)\Big\}\\
&:=&\tilde C.
\end{eqnarray*}
When $0\leq t\leq r_1<r_2<\tau$, similar to the case when $0\leq r_1<r_2\leq t<\tau$.
Therefore, ${\cal M}$ maps
 $C^0_{\tau,\alpha}((-\infty,+\infty),L^2(D,{\cal D}^{1,2}))$ to itself.$\hfill \hfill \sharp$.
 \vskip 0.5cm
  Define the set 
$$S:=C_\tau^0((-\infty,\infty),L^2(\Omega\times D))\cap L^\infty((-\infty,\infty),L^2(\Omega, H^1_0(D)))\cap C_{\tau,\alpha}^0((-\infty,\infty), L^2(D, {\cal D}^{1,2})).$$
 Define
$${\cal M}(S)|_{[0,\tau)}:= \{f|_{[0,\tau)}: f\in {\cal M}(S)\}.$$
\begin{lem}\label{lem2.3}
The set ${\cal M}(S)|_{[0,\tau)}$ is relatively compact in $C^0([0, \tau),L^2(\Omega\times D))$. 
\end{lem}
{\bf Proof:} With what we have proved in Lemma \ref{lem2.2}, we also need to prove that ${\cal D}_r {\cal M} (z)(t)$ is equicontinuous in $t$ in the space $L^2(D, {\cal D}^{1,2})$. We will consider several cases.\\
When $0\leq r\leq t_1<t_2<\tau$, for $z\in S$, 
\begin{eqnarray*}
&&\int_DE|{\cal D}_{r}{\cal M}(z)(t_2,x)-{\cal D}_{r}{\cal M}(z)(t_1,x)|^2dx\\
&\leq &C\int_D\Big\{E\big|\int_{-\infty}^{t_1}  \int_D\sum_{i=m+1}^\infty (e^{\mu_i(t_2-s)}-e^{\mu_i(t_1-s)})\phi_i(x)\phi_i(y)\nabla F^i(s, z(s)+Y_1(s))(y){\cal D}_{r} z(s,y)dyds\big|^2\\
&&\hskip 0.5cm+E\big|\int^{t_2}_{t_1} \sum_{i=m+1}^\infty e^{\mu_i(t_2-s)}\phi_i(x)\phi_i(y)\nabla F^i(s, z(s)+Y_1(s))(y){\cal D}_{r} z(s,y)dyds\big|^2\\
&&+E\big|\int^{\infty}_{t_2}\sum_{i=1}^m (e^{\mu_i(t_1-s)}-e^{\mu_i(t_2-s)})\phi_i(x)\phi_i(y)\nabla F^i(s, z(s)+Y_1(s))(y){\cal D}_{r} z(s,y)dyds\big|^2\\
&&\hskip 0.5cm+E\big|\int_{t_1}^{t_2} \sum_{i=1}^m e^{\mu_i(t_1-s)}\phi_i(x)\phi_i(y)\nabla F^i(s, z(s)+Y_1(s))(y){\cal D}_{r} z(s,y)dyds\big|^2\\
&&+E\big|\int^{t_1}_r \sum_{i=m+1}^\infty (e^{\mu_i(t_1-s)}-e^{\mu_i(t_2-s)})\phi_i(x)\phi_i(y)\nabla F^i(s, z(s)+Y_1(s))(y){\cal D}_{r} Y_1(s,y)dyds\big|^2\\
&&\hskip 0.5cm+E\big|\int_{t_1}^{t_2} \sum_{i=m+1}^\infty e^{\mu_i(t_2-s)}\phi_i(x)\phi_i(y)\nabla F^i(s, z(s)+Y_1(s))(y){\cal D}_{r} Y_1(s,y)dyds\big|^2\\
&&+E\big|\int^{r}_{-\infty} \sum_{i=m+1}^\infty (e^{\mu_i(t_2-s)}-e^{\mu_i(t_1-s)})\phi_i(x)\phi_i(y)\nabla F^i(s, z(s)+Y_1(s))(y){\cal D}_{r} Y_1(s,y)dyds\big|^2\\
&&+E\big|\int_{t_1}^{\infty}\sum_{i=1}^m (e^{\mu_i(t_1-s)}-e^{\mu_i(t_2-s)})\phi_i(x)\phi_i(y)\nabla F^i(s, z(s)+Y_1(s))(y){\cal D}_{r} Y_1(s,y)dyds\big|^2\\
&&\hskip 0.5cm+E\big|\int_{t_1}^{t_2} \sum_{i=1}^m e^{\mu_i(t_1-s)}\phi_i(x)\phi_i(y)\nabla F^i(s, z(s)+Y_1(s))(y){\cal D}_{r} Y_1(s,y)dyds\big|^2\Big\}dx\\
&:=&B_1+B_2+B_3+B_4+B_5+B_6+B_7+B_8+B_9.
\end{eqnarray*}
We will estimate them in the following steps. First, we have
\begin{eqnarray*}
&&B_1\\
&\leq&C\int_D E\Big|\int_{-\infty}^{t_1}\int_D \sum_{i=m+1}^\infty(e^{\mu_i(t_2-s)}-e^{\mu_i(t_1-s)})\phi_i(x)\phi_i(y)\nabla F^i(s,z(s)+Y_1(s))(y){\cal D}_r z(s,y)dyds\Big|^2dx\\
&\leq&C\sum_{i=m+1}^\infty E\int_{-\infty}^{t_1}\int_D |e^{\mu_i(t_2-s)}-e^{\mu_i(t_1-s)}|\cdot|\phi_i(y)|^2dyds\\
&&\hskip 1cm\cdot \int_{-\infty}^{t_1}\int_D |e^{\mu_i(t_2-s)}-e^{\mu_i(t_1-s)}| |\nabla F^i(s,z(s)+Y_1(s))(y)|^2|{\cal D}_r z(s,y)|^2dyds\\
&\leq & C|t_2-t_1| \cdot ||\nabla F||_\infty^2 \int_{-\infty}^{t_1}\int_De^{\mu_{m+1}(t_1-s)}E |{\cal D}_r z(s,y)|^2dyds\\
&\leq & C|t_2-t_1| \cdot ||\nabla F||_\infty^2 \Big[\int_{0}^{\tau}\int_D\sum_{i=0}^\infty e^{\mu_{m+1}(t_1-s+\tau+i\tau)}E |{\cal D}_r z(s,\theta_{-(i+1)\tau}\cdot,y)|^2dyds\\
&&\hskip 3.5cm+\int_{0}^{t_1}\int_De^{\mu_{m+1}(t_1-s)}E |{\cal D}_r z(s,y)|^2dyds\Big]\\
&\leq&-{C\over {\mu_{m+1}}}|t_2-t_1|\cdot ||\nabla F||_\infty^2 \alpha_1(\sum_{i=0}^\infty e^{\mu_{m+1}i\tau}+1).
\end{eqnarray*}
About $B_2$, we have 
\begin{eqnarray*}
&&B_2\\
&\leq& C\int_D E\Big|\int^{t_2}_{t_1}\int_D\sum_{i=m+1}^\infty (e^{\mu_i(t_2-s)}-e^{\mu_i(t_1-s)})\phi_i(x)\phi_i(y)\nabla F^i(s,z(s)+Y_1(s))(y){\cal D}_r z(s,y)dyds\Big|^2dx\\
&\leq&C\sum_{i=m+1}^\infty E\int^{t_2}_{t_1}\int_D |\phi_i(y)|^2dyds\cdot \int_{t_1}^{t_2}\int_D |\nabla F^i(s,z(s)+Y_1(s))(y)|^2|{\cal D}_r z(s,y)|^2dyds\\
&\leq&C|t_2-t_1|\cdot ||\nabla F||_\infty^2 \tau\alpha_1.
\end{eqnarray*}
Similar to $B_1$, we have
$$B_3\leq {C\over {\mu_{m}}}|t_2-t_1|\cdot ||\nabla F||_\infty^2 \alpha_1(\sum_{i=0}^\infty e^{-\mu_{m}i\tau}+1).$$
Similar to $B_2$, we have
$$B_4\leq C|t_2-t_1|\cdot ||\nabla F||_\infty^2 \tau\alpha_1.$$
About $B_5$,
\begin{eqnarray*}
B_5&\leq&  C\int_D E\Big|\int^{t_1}_r\int_D\sum_{i=m+1}^\infty (e^{\mu_i(t_2-s)}-e^{\mu_i(t_1-s)})\phi_i(x)\phi_i(y)\nabla F^i(s,z(s)+Y_1(s))(y)\\
&&\hskip 2.5cm \cdot \sum_{j=m+1}^\infty\phi_j(y)\sigma_j(r)dyds\Big|^2dx\\
&\leq& C\sum_{i=m+1}^\infty E\int_{r}^{t_1}\int_D |e^{\mu_i(t_2-s)}-e^{\mu_i(t_1-s)}|\cdot|\phi_i(y)|^2dyds\\
&&\cdot \int_{r}^{t_1}\int_D |e^{\mu_i(t_2-s)}-e^{\mu_i(t_1-s)}| |\nabla F^i(s,z(s)+Y_1(s))(y)|^2\Big|\sum_{j=m+1}^\infty \phi_j(y)\sigma_j(r)\Big|^2dyds\\
&\leq&C|t_2-t_1|\cdot ||\nabla F||_\infty^2\int_r^{t_1} \int_D\sum_{j=m+1}^\infty |\phi_j(y)|^2|\sigma_j(r)|^2dyds\\
&\leq&C|t_2-t_1|\cdot ||\nabla F||_\infty^2 \tau\sup_{s\in (-\infty,\infty)} \sum_{j=1}^\infty \sigma_j^2(s).
\end{eqnarray*}
About $B_6$,
\begin{eqnarray*}
B_6&\leq&C\int_D E\Big|\int^{t_2}_{t_1}\int_D\sum_{i=m+1}^\infty (e^{\mu_i(t_2-s)}-e^{\mu_i(t_1-s)})\phi_i(x)\phi_i(y) \nabla F^i(s,z(s)+Y_1(s))(y)\\
&&\hskip 2.5cm \cdot\sum_{j=m+1}^\infty\phi_j(y)\sigma_j(r)dyds\Big|^2dx\\
&\leq&C\sum_{i=m+1}^\infty E\int^{t_2}_{t_1}\int_D |\phi_i(y)|^2dyds\cdot \int_{t_1}^{t_2}\int_D |\nabla F^i(s,z(s)+Y_1(s))(y)|^2\Big|\sum_{j=m+1}^\infty\phi_j(y)\sigma_j(r)\Big|^2dyds\\
&\leq&C|t_2-t_1|\cdot ||\nabla F||_\infty^2\int_{t_1}^{t_2} \int_D\sum_{j=m+1}^\infty |\phi_j(y)|^2|\sigma_j(r)|^2dyds\\
&\leq&C|t_2-t_1|\cdot ||\nabla F||_\infty^2 \tau\sup_{s\in(-\infty,\infty)}  \sum_{j=1}^\infty|\sigma_j(s)|^2.
\end{eqnarray*}
Similarly, we have
\begin{eqnarray*}
B_7&\leq& C|t_2-t_1|\cdot ||\nabla F||_\infty^2(- {1\over {\mu_{m+1}}})\sup_{s\in (-\infty,\infty)} \sum_{j=1}^\infty \sigma_j^2(s),\\
B_8&\leq& C|t_2-t_1|\cdot ||\nabla F||_\infty^2 {1\over {\mu_{m}}}\sup_{s\in (-\infty,\infty)} \sum_{j=1}^\infty \sigma_j^2(s),\\
B_9&\leq& C|t_2-t_1|\cdot ||\nabla F||_\infty^2 \tau\sup_{s\in(-\infty,\infty)}  \sum_{j=1}^\infty|\sigma_j(s)|^2.
\end{eqnarray*}
Therefore, for any $z\in S$ and $0\leq r\leq t_1<t_2<\tau$,
\begin{eqnarray*}
\int_DE|{\cal D}_{r}{\cal M}(z)(t_2,x)-{\cal D}_{r}{\cal M}(z)(t_1,x)|^2dx\leq \tilde C |t_2-t_1|.
\end{eqnarray*}
When $0\leq t_1<r<t_2<\tau$, $z\in S$, similar as before, we can compute that
\begin{eqnarray*}
&&\int_DE|{\cal D}_{r}{\cal M}(z)(t_2,x)-{\cal D}_{r}{\cal M}(z)(t_1,x)|^2dx\\
&\leq &C\int_D\Big\{E\big|\int_{-\infty}^{t_2}\int_D \sum_{i=m+1}^\infty e^{\mu_i(t_2-s)}\phi_i(x)\phi_i(y)\nabla F^i(s, z(s)+Y_1(s))(y){\cal D}_{r} z(s,y)dyds\\
&&\hskip1cm-\int_{-\infty}^{t_1}\int_D \sum_{i=m+1}^\infty e^{\mu_i(t_1-s)}\phi_i(x)\phi_i(y)\nabla F^i(s, z(s)+Y_1(s))(y){\cal D}_{r} z(s,y)dyds\big|^2\\
&&+E\big|\int^{\infty}_{t_1} \int_D\sum_{i=1}^m e^{\mu_i(t_1-s)}\phi_i(x)\phi_i(y)\nabla F^i(s, z(s)+Y_1(s))(y){\cal D}_{r} z(s,y)dyds\\
&&\hskip 1cm-\int^{\infty}_{t_2}\int_D \sum_{i=1}^m e^{\mu_i(t_2-s)}\phi_i(x)\phi_i(y)\nabla F^i(s, z(s)+Y_1(s))(y){\cal D}_{r} z(s,y)dyds\big|^2\\
&&+E\big|\int^{t_2}_r \int_D\sum_{i=m+1}^\infty e^{\mu_i(t_2-s)}\phi_i(x)\phi_i(y)\nabla F^i(s, z(s)+Y_1(s))(y){\cal D}_{r} Y_1(s,y)dyds\big|^2\\
&&+E\big|\int_{t_1}^r \int_D \sum_{i=1}^m e^{\mu_i(t_1-s)}\phi_i(x)\phi_i(y)\nabla F^i(s, z(s)+Y_1(s))(y){\cal D}_{r} Y_1(s,y)dyds\big|^2\\
&&+E\big|-\int^{r}_{-\infty}\int_D  \sum_{i=m+1}^\infty e^{\mu_i(t_2-s)}\phi_i(x)\phi_i(y)\nabla F^i(s, z(s)+Y_1(s))(y){\cal D}_{r} Y_1(s,y)dyds\\
&&\hskip 1cm +\int^{t_1}_{-\infty}\int_D  \sum_{i=m+1}^\infty e^{\mu_i(t_1-s)}\phi_i(x)\phi_i(y)\nabla F^i(s, z(s)+Y_1(s))(y){\cal D}_{r} Y_1(s,y)dyds|^2\\
&&+E\big|-\int_{t_2}^{\infty} \int_D \sum_{i=1}^m e^{\mu_i(t_2-s)}\phi_i(x)\phi_i(y)\nabla F^i(s, z(s)+Y_1(s))(y){\cal D}_{r} Y_1(s,y)dyds\\
&&\hskip 1cm+\int_{r}^{\infty}\int_D  \sum_{i=1}^m e^{\mu_i(t_1-s)}\phi_i(x)\phi_i(y)\nabla F^i(s, z(s)+Y_1(s))(y){\cal D}_{r} Y_1(s,y)dyds\big|^2\Big\}dx\\
&\leq&C\Big\{(-{1\over {\mu_{m+1}}})(1+\sum_{i=0}^\infty e^{\mu_{m+1}i\tau})||\nabla F||_\infty^2\alpha_1(t_2-t_1)+{1\over {\mu_{m}}}(1+\sum_{i=0}^\infty e^{-\mu_{m}i\tau})||\nabla F||_\infty^2\alpha_1(t_2-t_1)\\
&&\hskip 0.5cm+2|t_2-t_1|\cdot ||\nabla F||_\infty^2 \tau\sup_{s\in(-\infty,\infty)}  \sum_{j=1}^\infty|\sigma_j(s)|^2\\
&&\hskip 0.5cm +(-{1\over {\mu_{m+1}}}+\tau)||\nabla F||_\infty^2\sum_{j=1}^\infty \sup_{s\in(-\infty,\infty)}  \sum_{j=1}^\infty|\sigma_j(s)|^2(t_2-t_1)\\
&&\hskip 0.5cm+({1\over {\mu_{m}}}+\tau)||\nabla F||_\infty^2\sum_{j=1}^\infty \sup_{s\in(-\infty,\infty)}  \sum_{j=1}^\infty|\sigma_j(s)|^2(t_2-t_1)\Big\}\\
&\leq & \hat C|t_2-t_1|.
\end{eqnarray*}
The case when $0\leq t_1<t_2<r<\tau$ is similar to the case when $0\leq r\leq t_1<t_2<\tau$.
Thus,
from the above arguments, by Theorem \ref{B-S}, ${\cal M}(S)|_{[0,\tau)}$ is relatively compact in  $C^0([0,\tau),L^2(\Omega\times D))$.$\hfill \hfill \sharp$
\vskip 0.5cm

From the periodicity of ${\cal M}(z)(t)$, we can prove
\begin{lem}
The set ${\cal M}(S)$ is relatively compact in $C_\tau^0((-\infty, +\infty), L^2(\Omega\times D))$.
\end{lem}
{\bf Proof:} From Lemma \ref{lem2.3}, we know for any sequence ${\cal M}(z_n)\in
S$, there exists a
subsequence, still denoted by ${\cal M}(z_n)$ and $Z^*\in
C^0([0,\tau), L^2(\Omega\times D))$ such that
\begin{eqnarray}\label{zhao1}
\sup\limits_{t\in [0,\tau)}\int_DE|{\cal M}(z_n)(t,\cdot,x)-Z^*(t,\cdot,x)|^2dx\to 0
\end{eqnarray}
as $n\to \infty$. Set for $\tau\leq t<2\tau$,
$$Z^*(t,\omega,x)=Z^*(t-\tau,\theta_{\tau}\omega,x).$$
Noting
\begin{eqnarray*}
{\cal M}(z_n)(t,\theta_{\tau}\omega,x)={\cal M}(z_n)(t+\tau,\omega,x),
\end{eqnarray*}
from (\ref{zhao1}), and the probability preserving property of $\theta$, we have
\begin{eqnarray*}
\sup\limits_{t\in [\tau,2\tau)}\int_DE|{\cal M}(z_n)(t,\cdot,x)-Z^*(t,\cdot,x)|^2dx
&=&\sup\limits_{t\in [0,\tau)}\int_DE|{\cal M}(z_n)(t+\tau ,\cdot,x)-Z^*(t+\tau,\cdot,x)|^2dx\\
&=&\sup\limits_{t\in [0,\tau)}\int_DE|{\cal M}(z_n)(t,\theta
_{\tau}\cdot,x)-Z^*(t,\theta_{\tau}\cdot,x)|^2dx\\
&=&\sup\limits_{t\in [0,\tau)}\int_DE|{\cal M}(z_n)(t,\cdot,x)-Z^*(t,\cdot,x)|^2dx\\
&\to& 0.
\end{eqnarray*}
Similarly one can prove that
\begin{eqnarray}
&&\sup\limits_{t\in [0,\tau)}\int_DE|{\cal M}(z_n)(t+m\tau ,\cdot,x)-Z^*(t+m\tau,\cdot,x)|^2dx\\
&=&\sup\limits_{t\in [0,\tau)}\int_DE|{\cal M}(z_n)(t,\cdot,x)-Z^*(t,\cdot,x)|^2dx\to 0,
\end{eqnarray}
for any $m\in \{0, \pm1, \pm2,\cdots\}$. Therefore
\begin{eqnarray*}
\sup\limits_{t\in (-\infty,+\infty)}\int_DE|{\cal M}(z_n)(t,
\cdot,x)-Z^*(t,\cdot,x)|^2dx\to 0,
\end{eqnarray*}
as $n\to \infty$. Therefore ${\cal M}(S)$ is relatively compact in $C_\tau^0((-\infty, +\infty), L^2(\Omega\times D))$.$\hfill \hfill \sharp$
\\

{\bf Proof of Theorem \ref{aug20b}:}
From the above four lemmas, according to the generalized
Schauder's fixed point theorem, ${\cal M}$ has a fixed point in
$C_{\tau}^0((-\infty, +\infty), L^2(\Omega\times D))$. That is to say
there exists a solution $Z\in
C_{\tau}^0((-\infty, +\infty), L^2(\Omega\times D))$ of equation
(\ref{zhao5}) such that for any $t\in (-\infty,+\infty)$,
$Z(t+\tau,\omega,x)=Z(t,\theta_{\tau}\omega,x)$. Then $Y=Z+Y_1$ is the
desired solution of
(\ref{sep17a}).  Moreover, $Y(t+\tau,\omega,x)=Y(t,\theta_{\tau}\omega,x).$\hfill \hfill $\sharp$
\\

Now we consider the semilinear
stochastic differential equations with the additive noise of the form
\begin{eqnarray}\label{(bo50)}
du(t,x)&=&[{\cal L}u(t,x)+F(u(t,x))]dt+\sum_{k=1}^\infty \sigma_k\phi_k(x)W^k(t),\\
u(0)&=&\psi\in L^2(D), \nonumber\\
u(t)|_{\partial D}&=&0\nonumber,
\end{eqnarray}
for $t\geq 0$. Here $F$ and $\sigma_k$ do not depend on time $t$, that is to say, $\tau$ in Condition (P) can be chosen as an arbitrary real number. We have a similar variation of constant representation
to (\ref{sep17a}). The difference is that for this equation, we have a cocycle.
Similar to Theorem \ref{aug20d}, we can prove the following theorem. But we do not give the proof here.

\begin{thm}\label{09-27a}
Assume Cauchy problem (\ref{(bo50)}) has a unique
solution $u(t,\omega, x)$ and the coupled forward-backward infinite horizon
stochastic integral equation
\begin{eqnarray}\label{(bo1)}
Y(\omega)&=&\int^{0}_{-\infty}T_{-s}P^-F(Y(\theta_s\omega))ds-\int^{\infty}_{0}T_{-s}P^+F(Y(\theta_s
\omega))ds\nonumber\\
&&+(\omega)\sum_{k=1}^\infty\int^{0}_{-\infty} \sigma_kT_{-s}P^-\phi_kW^k(s)-(\omega)\sum_{k=1}^\infty\int^{\infty}_{0} \sigma_kT_{-s}P^+\phi_kW^k(s)
\end{eqnarray}
has one solution $Y: \Omega\rightarrow  L^2(D)$, then $Y$ is a stationary solution of
equation (\ref{(bo50)}) i.e.
\begin{eqnarray}
u(t, Y(\omega),\omega)=Y(\theta_{t} \omega) \ \ {\rm for \ any} \ \
t\geq 0
\ \ \ \ a.s.
\end{eqnarray}
Conversely, if equation (\ref{(bo50)}) has a stationary solution $Y: \Omega\rightarrow L^2(D)$
which is tempered from above,
then $Y$ is a solution of the coupled forward-backward infinite horizon
stochastic integral equation (\ref{(bo1)}).
\end{thm}
\begin{thm}\label{aug20a}
Assume the same onditions on ${\cal L}$ as in Theorem \ref {aug20b} and $\sum_{k=1}^{\infty} \sigma_k^2<\infty$. Let $F:  R\to R$ be a continuous map, globally bounded and $\nabla F$ being globally bounded.
 Then there exists at least one $\mathcal{F}$-measurable map
$Y:\Omega\rightarrow  L^2(D)$ satisfying (\ref{(bo1)}).
\end{thm}
{\bf Proof:} Set the $\mathcal{F}$-measurable map
$Y_{1}:\Omega\rightarrow  L^2(D)$
\begin{eqnarray}\label{(bo2)}
Y_{1}(\omega)=(\omega)\sum_{k=1}^\infty\int^{0}_{-\infty} \sigma_kT_{-s}P^-\phi_kW^k(s)-(\omega)\sum_{k=1}^\infty\int^{\infty}_{0} \sigma_kT_{-s}P^+\phi_kW^k(s).
\end{eqnarray}
Then we have
\begin{eqnarray*}
Y_{1}(\theta_t\omega)&=&(\theta_t\omega)\sum_{k=1}^\infty\int^{0}_{-\infty} \sigma_kT_{-s}P^-\phi_kW^k(s)-(\theta_t\omega)\sum_{k=1}^\infty\int^{\infty}_{0} \sigma_kT_{-s}P^+\phi_kW^k(s)\\
&=&(\omega)\sum_{k=1}^\infty\int^{t}_{-\infty} \sigma_kT_{-s}P^-\phi_kW^k(s)-(\omega)\sum_{k=1}^\infty\int^{\infty}_{t} \sigma_kT_{-s}P^+\phi_kW^k(s).
\end{eqnarray*}
We need to solve the equation
\begin{eqnarray}
&&Z(t,\omega)\nonumber\\
&=&\int^{t}_{-\infty}T_{t-s}P^-F(Z(s,\omega))+Y_1(\theta_s\omega))ds-\int^{\infty}_{t}T_{t-s}P^+F(Z(s,\omega)+Y_1(\theta_s\omega)))ds.
\end{eqnarray}
For this, define
\begin{eqnarray*}
&&C_s^0((-\infty, +\infty), L^2(\Omega\times D))\nonumber\\
 &:=&\{f\in C^0((-\infty, +\infty), L^2(\Omega\times D)):
 \ {\rm for \ any} \ \
t\in (-\infty,\infty),\
f(t,\omega,x)=f(0,\theta_t\omega,x) \},\nonumber
\end{eqnarray*}
We now define for any $z\in
C_s^0((-\infty, +\infty), L^2(\Omega\times D))$,
\begin{eqnarray}\label{2.20}
{\cal M}(z)(t,\omega)&=&\int^{t}_{-\infty}T_{t-s}P^-F(z(s,\omega)+Y_{1}(\theta_s\omega))ds\nonumber
\\
&&-\int^{+\infty}_{t}T_{t-s}P^+F(z(s,\omega)+Y_{1}(\theta_{s}\omega))ds.
\end{eqnarray}
It's easy to see that
\begin{eqnarray*}
&&{\cal M}(z)(0,\theta _t\omega)\\
&=&\int^{0}_{-\infty}T_{-s}P^-F(z(s,\theta _t\omega)+Y_{1}(\theta_{s+t}\omega))ds-\int^{+\infty}_{0}T_{-s}P^+F(z(s,\theta
_t\omega)+Y_{1}(\theta_{s+t}\omega))ds\\
&=&\int^{0}_{-\infty}T_{-s}P^-F(z(s+t,\omega)+Y_{1}(\theta_{s+t}\omega))ds-\int^{+\infty}_{0}T_{-s}P^+F(z(s+t,\omega)+Y_{1}(\theta_{s+t}\omega,x))ds\\
&=&\int^{t}_{-\infty}T_{t-s}P^-F(z(s,\omega)+Y_{1}(\theta_{s}\omega))ds-\int^{+\infty}_{t}T_{t-s}P^+F(z(s,\omega)+Y_{1}(\theta_{s}\omega))ds\\
&=& {\cal M}(z)(t,\omega).
\end{eqnarray*}
By the similar method in the proof of Lemma \ref{lem2.1}, we can see that the ${\cal M}$ defined in (\ref{2.20})  maps $C_s^0((-\infty, +\infty), L^2(\Omega\times D))\to C_s^0((-\infty, +\infty), L^2(\Omega\times D))$ is a continuous map. Moreover ${\cal M}$ maps $ C_s^0((-\infty, +\infty), L^2(\Omega\times D))$ into $ C_s^0((-\infty, +\infty), L^2(\Omega\times D))\cap L^\infty((-\infty, +\infty), L^2(\Omega,H^1_0(D))).$ For a fixed $T>0$, define
\begin{eqnarray*}
&&C^0_{T,\alpha}((-\infty,+\infty),L^2(D,{\cal D}^{1,2})) \\
&:=&\{ f\in C_{T}^0((-\infty,+\infty), L^2(\Omega\times D)):\ f|_{[0,T)}\in C^0([0,T),L^2(D,{\cal D}^{1,2})),\\
&& i.e.\  ||f||^2=\sup_{t\in [0,T)}\int_D||f(t,x)||^2_{1,2}dx<\infty, {\rm \ and\  for \ any}\  t,r\in [0,T),\ i=0,\pm 1, \pm 2, \cdots \nonumber\\
&&\int_DE|{\cal D}_r f(t,\theta_{iT}\cdot, x)|^2dx\leq \alpha_r(t), \sup_{s,r_1, r_2 \in [0,T)}{{\int_DE|{\cal D}_{r_1} f(s,\theta_{iT}\cdot,x)-{\cal D}_{r_2} f(s,\theta_{iT}\cdot,x)|^2dx}\over {|r_1-r_2|}}<\infty\}.
\end{eqnarray*}
Here $\alpha_r(t)$ is the solution of integral equation (see page 324 in \cite{polyanin}) 
\begin{eqnarray}\label{eqn2.15}
\alpha_r(t)=A\int_{r-2T}^{r+2T} e^{-\beta |t-s|}\alpha_r(s)ds +B,
\end{eqnarray}
where
\begin{eqnarray*}
&&A=C||\nabla F||^2_\infty(-{1\over {\mu_{m+1}}}\sum_{i=0}^\infty e^{\mu_{m+1}iT}+{1\over {\mu_{m}}}\sum_{i=0}^\infty e^{-\mu_{m}iT}), \\
&&B=C||\nabla F||^2_\infty \sup_{s\in (-\infty,\infty)} \sum_{j=1}^\infty \sigma_j^2(s)({1\over {\mu^2_{m+1}}}+{1\over {\mu^2_{m}}}),\  \beta=min\{-\mu_{m+1},\mu_m\}. 
\end{eqnarray*}
And similar to Lemma \ref{lem2.2} we can get ${\cal M}$ maps $C^0_{T,\alpha}((-\infty,+\infty),L^2(D,{\cal D}^{1,2}))$ into itself.  
 Define the set 
$$S:=C_T^0((-\infty,\infty),L^2(\Omega\times D))\cap L^\infty((-\infty,\infty),L^2(\Omega, H^1_0(D)))\cap C_{T,\alpha}^0((-\infty,\infty),L^2(D,{\cal D}^{1,2})).$$
 Similar to Lemma \ref{lem2.3} we can prove the set ${\cal M}(S)|_{[0,T)}$ is relatively compact in $C^0([0, T),L^2(\Omega\times D))$. 
We need to prove that ${\cal M}( S)$ is relatively compact in $C_s^0((-\infty, +\infty), L^2(\Omega\times D))$.
 Note also for any sequence ${\cal M}(z_n)\in{\cal M}(S)$, there exists a
subsequence, still denoted by ${\cal M}(z_n)$ and $Z^*\in
C^0([0,T), L^2(\Omega\times D))$ such that
$$\int_DE|{\cal M}(z_n)(0,\cdot,x)-Z^*(\cdot,x)|^2dx\to 0, \ \ {\rm as} \ \ n
\to \infty.$$
Define
\begin{eqnarray*}
Z^*(t,\omega,x)=Z^*(0,\theta_t\omega,x).
\end{eqnarray*}
Noting
\begin{eqnarray*}
{\cal M}(z_n)(0,\theta_{t}\omega,x)={\cal M}(z_n)(t,\omega,x),
\end{eqnarray*}
and by the probability preserving property of $\theta$, we have
\begin{eqnarray*}
\sup_{t\in (-\infty,\infty)}\int_DE|{\cal M}(z_n)(t,\cdot,x)-Z^*(\theta_t\cdot,x)|^2dx
&=&\sup_{t\in (-\infty,\infty)}\int_DE|{\cal M}(z_n)(0,\theta_t\cdot)-Z^*(\theta_t\cdot)|^2dx\\
&=&\int_DE|{\cal M}(z_n)(0,\cdot,x)-Z^*(\cdot,x)|^2dx\\
&\to& 0, \ \ {\rm as} \ \ n
\to \infty.
\end{eqnarray*}
 So
${\cal M}(S)$ is relatively compact in $C_s^0((-\infty, +\infty), L^2(\Omega\times D))$.
Therefore, according to generalized Schauder's fixed point theorem, ${\cal M}$ has a
fixed point in $C_s^0((-\infty, +\infty), L^2(\Omega\times D))$. That is to
say that there exists $Z\in C_s^0((-\infty, +\infty), L^2(\Omega\times D))$
such that for any $t\in (-\infty,+\infty)$,
$Z(t,\omega)=Z(0,\theta_t\omega)$ and
\begin{eqnarray*}
Z(0,\theta_t\omega)&=&\int^{t}_{-\infty}T_{t-s}P^-F(Z(0,\theta_{s}\omega)+Y_{1}(\theta_{s}\omega))ds-\int^{+\infty}_{t}T_{t-s}P^+F(Z(0,\theta_{s}\omega)+Y_{1}(\theta_{s}\omega))ds.
\end{eqnarray*}
Finally, we add $Y_1$ defined by the integral equation (\ref{(bo2)})
to the above equation and also assume
$$Y(\omega):=Z(0,\omega)+Y_1(\omega).$$
It's easy to see that $Y(\omega,x)$ satisfies (\ref{(bo1)}). \hfill $\sharp$\\
\\
{\bf Acknowledgements}. 
We would like to thank the referee for very useful comments and pointing out to us the references \cite{cho1,cho2}, and \cite{klunger}.

\end{document}